\pdfoutput=1 
\documentclass[11pt]{amsart}
\usepackage{amssymb,amsthm,enumitem,colonequals,tikz-cd,microtype,rotating}
\usepackage[normalem]{ulem} 

\usepackage[osf]{Baskervaldx}
\usepackage[baskervaldx]{newtxmath}
\usepackage[cal=boondoxo]{mathalfa} 

\usepackage[top=3.75cm, bottom=3cm, left=3.5cm, right=3.5cm]{geometry}

\usepackage{xcolor}
\colorlet{darkblue}{blue!55!black}
\colorlet{darkcyan}{cyan!50!black}
\colorlet{darkgreen}{green!60!black}

\PassOptionsToPackage{hyphens}{url}
\usepackage{hyperref}
\hypersetup{
  colorlinks=true,
  linkcolor=darkblue,
  urlcolor=darkcyan,
  citecolor=darkgreen,
}

\def\eqref#1{\textcolor{darkblue}{(\ref{#1})}}

\usepackage[nameinlink]{cleveref} 
\Crefformat{section}{#2\S#1#3}
\Crefmultiformat{section}{#2\S\S#1#3}{ and~#2#1#3}{, #2#1#3}{, and~#2#1#3}
\Crefname{equation}{Diagram}{Diagrams}

\usepackage[pagewise]{lineno}
\let\oldequation\equation
\let\oldendequation\endequation
\renewenvironment{equation}{\linenomathNonumbers\oldequation}{\oldendequation\endlinenomath}
\expandafter\let\expandafter\oldequationstar\csname equation*\endcsname
\expandafter\let\expandafter\oldendequationstar\csname endequation*\endcsname
\renewenvironment{equation*}{\linenomathNonumbers\oldequationstar}{\oldendequationstar\endlinenomath}
\let\oldalign\align
\let\oldendalign\endalign
\renewenvironment{align}{\linenomathNonumbers\oldalign}{\oldendalign\endlinenomath}
\expandafter\let\expandafter\oldalignstar\csname align*\endcsname
\expandafter\let\expandafter\oldendalignstar\csname endalign*\endcsname
\renewenvironment{align*}{\linenomathNonumbers\oldalignstar}{\oldendalignstar\endlinenomath}

\theoremstyle{plain}
\newtheorem{theorem}{Theorem}[section]
\newtheorem{lemma}[theorem]{Lemma}
\newtheorem{corollary}[theorem]{Corollary}
\newtheorem{proposition}[theorem]{Proposition}

\theoremstyle{definition}
\newtheorem{definition}[theorem]{Definition}
\newtheorem{example}[theorem]{Example}
\newtheorem{remark}[theorem]{Remark}
\newtheorem{setup}[theorem]{Setup}
\newtheorem{question}[theorem]{Question}

\newtheorem*{ack}{Acknowledgments}

\AddToHook{env/conjecture/begin}{\crefalias{theorem}{conjecture}}
\AddToHook{env/lemma/begin}{\crefalias{theorem}{lemma}}
\AddToHook{env/corollary/begin}{\crefalias{theorem}{corollary}}
\AddToHook{env/proposition/begin}{\crefalias{theorem}{proposition}}
\AddToHook{env/definition/begin}{\crefalias{theorem}{definition}}
\AddToHook{env/remark/begin}{\crefalias{theorem}{remark}}
\AddToHook{env/setup/begin}{\crefalias{theorem}{setup}}
\AddToHook{env/example/begin}{\crefalias{theorem}{example}}
\AddToHook{env/question/begin}{\crefalias{theorem}{question}}

\numberwithin{equation}{section}
\numberwithin{theorem}{section}

\title[Fiberwise criteria for Fourier--Mukai equivalences]{Fiberwise criteria for Fourier--Mukai equivalences}

\author[E.~Guisado Villalgordo]{El\'{i}as Guisado Villalgordo}
\address{E.~Guisado Villalgordo,
BCAM -- Basque Center for Applied Mathematics, Mazarredo 14, 48009
Bilbao, Basque Country -- Spain}
\email{eguisado@bcamath.org}

\author[P.~Lank]{Pat Lank}
\address{P.~Lank,
Dipartimento di Matematica `F.~Enriques', Universit\`{a} degli Studi di Milano, Via Cesare Saldini 50, 20133 Milano, Italy}
\email{lankp.math@gmail.com}

\author[K.~Manali Rahul]{Kabeer Manali Rahul}
\address{K.~Manali Rahul,
Faculty of Mathematics and Physics, Charles University, Sokolovsk\'{a} 83, 186 75 Praha 8, Czech Republic}
\email{kabeermr.maths@gmail.com}

\author[N.~Pavic]{Nebojsa Pavic}
\address{N.~Pavic,
Institut f\"{u}r Mathematik und Wissenschaftliches Rechnen,
Universit\"{a}t Graz,
Graz, Austria}
\email{nebojsa.pavic@uni-graz.at}

\date{\today}



\begin{document}

\begin{abstract}
	We prove fiberwise criteria for full faithfulness and equivalences of integral transforms between proper flat families of schemes over a separated Noetherian base.
	The kernels are bounded pseudocoherent complexes which are relatively perfect over each factor, and need not be perfect.
	A key ingredient is the calculus of mates, which controls the units and counits of adjoint integral transforms under base change.
	In particular, we give a complete argument for a fiberwise criterion of Hern\'{a}ndez Ruip\'{e}rez--L\'{o}pez Mart\'{i}n--Sancho de Salas, and extend it to arbitrary separated Noetherian bases.
	As applications, we extend a result of Orlov on field extensions to singular varieties, and show that over a Noetherian local ring these properties can be detected on the special fiber.
\end{abstract}

\maketitle

\tableofcontents

\section{Introduction}
\label{sec:intro}

\subsection{Motivation}
\label{sec:intro_motivation}

This article addresses the following question.

\begin{question}
	Can full faithfulness and equivalences of integral transforms between families be detected on fibers?
\end{question}

We make this precise.
Suppose $f_i\colon Y_i \to S$ are proper flat morphisms of Noetherian schemes.
Consider the fibered square
\begin{displaymath}
	\begin{tikzcd}
		{Y_1\times_S Y_2} & {Y_2} \\
		{Y_1} & {S.}
		\arrow["{f_1^\prime}", from=1-1, to=1-2]
		\arrow["{f_2^\prime}"', from=1-1, to=2-1]
		\arrow["{f_2}", from=1-2, to=2-2]
		\arrow["{f_1}"', from=2-1, to=2-2]
	\end{tikzcd}
\end{displaymath}
For any $K\in D_{\operatorname{qc}}(Y_1\times_S Y_2)$, the \textbf{integral transform with kernel} $K$ is the functor $\Phi_K\colon D_{\operatorname{qc}}(Y_1) \to D_{\operatorname{qc}}(Y_2)$ given by $A\mapsto \mathbf{R}(f_1^\prime)_\ast(\mathbf{L}(f_2^\prime)^\ast A\otimes^{\mathbf{L}} K)$.
We say $Y_1$ and $Y_2$ are \textbf{Fourier--Mukai $S$-partners} if there exists $K\in D^b_{\operatorname{coh}}(Y_1\times_S Y_2)$ such that $\Phi_K$ restricts to an equivalence $D^b_{\operatorname{coh}}(Y_1)\to D^b_{\operatorname{coh}}(Y_2)$.

Let $p\in S$.
Denote by $t\colon \operatorname{Spec}(\kappa(p))\to S$ the canonical morphism and by $t^\prime \colon Y_1\times_S Y_2\times_S \operatorname{Spec}(\kappa(p)) \to Y_1\times_S Y_2$ its base change.
Then the \textbf{fiber} of $\Phi_K$ at $p$ is the integral transform
\begin{displaymath}
	\Phi_{\mathbf{L}(t^\prime)^\ast K}\colon D_{\operatorname{qc}}(Y_1\times_S \operatorname{Spec}(\kappa(p))) \to D_{\operatorname{qc}}(Y_2\times_S \operatorname{Spec}(\kappa(p))).
\end{displaymath}
The question above asks whether $\Phi_K$ is fully faithful (resp.\ an equivalence) if its fibers are.

\subsection{Background}
\label{sec:intro_background}

Orlov proved that an integral transform which induces an equivalence between smooth projective varieties remains an equivalence after base change along any field extension \cite[Lemma 2.12]{Orlov:2002}.
The argument uses that the units and counits of an adjoint pair of integral transforms with perfect kernels are induced by morphisms between kernels \cite{Anno/Logvinenko:2012,Orlov:2002}.
This relies on the perfectness of the kernels.
See \cite[Theorem 3.1 and Appendix A, Eq.\ (A.18)]{Anno/Logvinenko:2012}.
It does not apply directly to singular varieties, where bounded coherent complexes need not be perfect.

For possibly singular varieties over an algebraically closed field, criteria determining when an integral transform is fully faithful or an equivalence were established in \cite{Ruiperez/Martin/deSalas:2007, Ruiperez/Hernandez/Martin/SanchodeSalas:2009}.
A natural problem is to extend these results to arbitrary fields, as well as to arithmetic situations (e.g.\ projective schemes over a DVR of mixed characteristic).
The present work grew out of the fiberwise criterion \cite[Proposition 2.15]{Ruiperez/Hernandez/Martin/SanchodeSalas:2009}, whose proof leaves a compatibility between units, counits, and base change unverified.
See \Cref{sec:intro_complete}.

\subsection{Results}
\label{sec:intro_results}

We work with bounded pseudocoherent kernels which are relatively perfect over each factor (see \Cref{sec:prelim_integral_transforms}).
These kernels need not be perfect.
By \cite[Propositions 3.2 and 3.6]{Dutta/Lank/ManaliRahul:2025}, relative perfectness over $Y_1$ (resp.\ $Y_2$) is equivalent to the integral transform preserving bounded pseudocoherent (resp.\ perfect) complexes.
In particular, for such kernels, $\Phi_K$ restricts to functors on $D^b_{\operatorname{coh}}$ and on $\operatorname{Perf}$, and we say that $\Phi_K$ is fully faithful (resp.\ an equivalence) on $D^b_{\operatorname{coh}}$ if its restriction is.

We also record that, for a proper morphism $f\colon Y\to X$ of Noetherian schemes, a pseudocoherent complex on $Y$ is $f$-perfect (i.e.\ relatively perfect over $X$) if, and only if, it has finite tor-dimension as an object of $D(f^{-1}\mathcal{O}_X)$.
See \Cref{cor:relative_perfectness_tor_dimension}.
Thus, for proper flat morphisms, our notion of relative perfectness agrees with the notion of relatively perfect objects in \cite[\href{https://stacks.math.columbia.edu/tag/0DI0}{Tag 0DI0}]{StacksProject}.

\subsubsection{Fiberwise criteria}
\label{sec:intro_fiberwise}

Our main result addresses the question above.

\begin{theorem}
	\label{introthm:fiberwise}
	Let $S$ be a separated Noetherian scheme.
	Consider proper flat morphisms $f_1\colon Y_1 \to S$ and $f_2 \colon Y_2 \to S$.
	Suppose $K\in D^b_{\operatorname{coh}}(Y_1\times_S Y_2)$ is relatively perfect over each $Y_i$.
	Then $\Phi_K$ is fully faithful (resp.\ an equivalence) on $D^b_{\operatorname{coh}}$ if, and only if, for every closed point $p\in S$ with associated closed immersion $t\colon \operatorname{Spec}(\kappa(p))\to S$, the integral transform $\Phi_{\mathbf{L}(t^\prime)^\ast K}$ is fully faithful (resp.\ an equivalence) on $D^b_{\operatorname{coh}}$, where $t^\prime\colon Y_1\times_S Y_2\times_S \operatorname{Spec}(\kappa(p))\to Y_1\times_S Y_2$ is the base change of $t$.
\end{theorem}

See \Cref{thm:fiberwise}, which also shows that full faithfulness and equivalences are preserved by affine base change, and that they can be detected after base change along any affine surjective morphism.
No smoothness assumptions are made on the $Y_i$ or on their fibers.
The separatedness of $S$ is only used to describe the adjoints of integral transforms via \cite[Lemmas 3.7 and 3.10]{Ballard:2009}.
See \Cref{sec:base_change_adjoints}.

\subsubsection{Completeness}
\label{sec:intro_complete}

For schemes of finite type over an algebraically closed field with $f_1$ locally projective, the statement of \Cref{introthm:fiberwise} appears as \cite[Proposition 2.15]{Ruiperez/Hernandez/Martin/SanchodeSalas:2009}.
The argument in loc.\ cit.\ uses that derived pullback and pushforward along the fibers carry the unit and counit of the adjunction for $\Phi_K$ to those of the fiberwise adjunction.
However, this compatibility is not justified there.
It is not formal, as the mate of a natural isomorphism need not be an isomorphism (see \Cref{rem:mates_of_isomorphisms}).

The main technical step of the present work is establishing this compatibility.
See \Cref{prop:left_mate_isomorphism} and \Cref{rem:unit_counit_compatibility}.
In particular, \Cref{introthm:fiberwise} provides a complete proof of \cite[Proposition 2.15]{Ruiperez/Hernandez/Martin/SanchodeSalas:2009}, and extends it to proper (not necessarily projective) flat morphisms over an arbitrary separated Noetherian base.
In the setting of algebraic stacks and perfect kernels, full faithfulness can be detected on geometric fibers over closed points \cite[Proposition 5.7]{Hall/Priver:2024}.
\Cref{introthm:fiberwise} allows for nonperfect kernels, and it treats equivalences.

\subsubsection{Applications}
\label{sec:intro_applications}

The ascent and descent statements of \Cref{thm:fiberwise} also extend Orlov's result \cite[Lemma 2.12]{Orlov:2002} to the singular setting, and our proof is independent of loc.\ cit.

\begin{corollary}
	\label{introcor:orlov}
	Let $Y_1$ and $Y_2$ be proper schemes over a field $k$.
	Suppose $K\in D^b_{\operatorname{coh}}(Y_1\times_k Y_2)$ is relatively perfect over each $Y_i$.
	Let $L/k$ be a field extension and $t^\prime\colon Y_1\times_k Y_2\times_k \operatorname{Spec}(L) \to Y_1\times_k Y_2$ the canonical morphism.
	Then $\Phi_K$ is fully faithful (resp.\ an equivalence) on $D^b_{\operatorname{coh}}$ if, and only if, $\Phi_{\mathbf{L}(t^\prime)^\ast K}$ is fully faithful (resp.\ an equivalence) on $D^b_{\operatorname{coh}}$.
\end{corollary}

See \Cref{cor:orlov}.
Combined with \Cref{introthm:fiberwise}, it shows that full faithfulness and equivalences can be detected on the geometric fibers over closed points of $S$, i.e.\ after base change along $\operatorname{Spec}(\overline{\kappa(p)})\to S$ for every closed point $p\in S$.
Over a local base, \Cref{introthm:fiberwise,thm:fiberwise} take the following form.

\begin{example}
	\label{ex:local_ring}
	Let $(R,\mathfrak{m},k)$ be a Noetherian local ring.
	Suppose $f_i\colon Y_i\to \operatorname{Spec}(R)$ are proper flat morphisms and $K\in D^b_{\operatorname{coh}}(Y_1\times_R Y_2)$ is relatively perfect over each $Y_i$.
	Let $t\colon \operatorname{Spec}(k)\to \operatorname{Spec}(R)$ and $s\colon \operatorname{Spec}(\widehat{R}) \to \operatorname{Spec}(R)$ be the morphisms induced respectively by the quotient homomorphism $R\to k$ and the $\mathfrak{m}$-adic completion $R\to\widehat{R}$.
	Denote by $t^\prime$ and $s^\prime$ their base changes to $Y_1\times_R Y_2$.
	By \Cref{introthm:fiberwise} and, as $s$ is affine and faithfully flat, \Cref{thm:fiberwise}, the following are equivalent:
	\begin{itemize}
		\item $\Phi_K$ is fully faithful (resp.\ an equivalence) on $D^b_{\operatorname{coh}}$
		\item $\Phi_{\mathbf{L}(t^\prime)^\ast K}$ is fully faithful (resp.\ an equivalence) on $D^b_{\operatorname{coh}}$
		\item $\Phi_{\mathbf{L}(s^\prime)^\ast K}$ is fully faithful (resp.\ an equivalence) on $D^b_{\operatorname{coh}}$.
	\end{itemize}
	In other words, full faithfulness and equivalences can be detected along the special fiber, or after passage to the completion.
\end{example}

\subsubsection{Methods}
\label{sec:intro_methods}

Our approach uses the calculus of mates to establish compatibility of adjunctions with base change, and does not require units and counits to be induced by morphisms between kernels.
The calculus of mates produces canonical comparison morphisms, but does not by itself imply that they are isomorphisms.
Using Grothendieck duality and relative perfectness, we show that the relevant comparison morphisms are isomorphisms.
See \Cref{lem:base_change_relative_dual} and \Cref{prop:left_mate_isomorphism}.
Mates are also used by Bergh--Schn\"{u}rer \cite{Bergh/Schnurer:2020} to study semiorthogonal decompositions under base change.
In loc.\ cit., perfectness of the kernels is used, see e.g.\ the proof of \cite[Proposition 4.9]{Bergh/Schnurer:2020}.

\begin{ack}
	The authors thank Rom\'{a}n Ahumada Gialanella, Matthew R.\ Ballard, Arend Bayer, Ana Cristina L\'{o}pez Mart\'{i}n, Jack Hall, Andres Fernandez Herrero, Emilio Franco, Ian Grojnowski, Fei Peng, Fernando Sancho de Salas, Dylan Spence, and Paolo Stellari for discussions and suggestions.
	Lank was supported by the ERC Advanced Grant 101095900-TriCatApp.
	Guisado Villalgordo was supported under the grant PRE2022-104796 by the Spanish Ministry of Science and Innovation.
	Manali Rahul thanks the Max Planck Institute for Mathematics in Bonn for hospitality and financial support.
\end{ack}

\section{Preliminaries}
\label{sec:prelim}

By a `natural' morphism, we mean something which is functorial.
By a `canonical' morphism, we mean something which is constructed and might not necessarily be functorial.

\subsection{Generation}
\label{sec:prelim_generation}

Let $\mathcal{T}$ be a triangulated category with shift functor $[1]\colon \mathcal{T} \to \mathcal{T}$.
Fix a subcategory $\mathcal{S}\subseteq \mathcal{T}$.
Denote by $\operatorname{add}(\mathcal{S})$ the strictly full subcategory of $\mathcal{T}$ consisting of all direct summands of objects of the form $\oplus_{n\in \mathbb{Z}} S_n^{\oplus r_n}[n]$ where $S_n \in \mathcal{S}$ for all $n$ and finitely many $r_n$ are nonzero.
Inductively, let $\langle \mathcal{S}\rangle_0$ consist of all zero objects in $\mathcal{T}$, $\langle \mathcal{S}\rangle_1 := \operatorname{add}(\mathcal{S})$, and
\begin{displaymath}
	\langle \mathcal{S}\rangle_n := \operatorname{add}(\{ \operatorname{cone}(\phi) : \phi \in \operatorname{Hom}(\langle \mathcal{S}  \rangle_{n-1}, \langle \mathcal{S}  \rangle_1) \}) \textrm{ for } n>1.
\end{displaymath}
Set $\langle\mathcal{S} \rangle:=\cup_{n\geq 0} \langle\mathcal{S} \rangle_n$.
It may be checked that $\langle\mathcal{S} \rangle$ is the smallest \textbf{thick} (that is, triangulated and closed under direct summands) subcategory of $\mathcal{T}$ containing $\mathcal{S}$.
If $\mathcal{S}$ consists of a single object $G$, we write $\langle G \rangle_n$ and $\langle G \rangle$.
An object $G\in \mathcal{T}$ is called a \textbf{classical generator} for $\mathcal{T}$ if $\langle G \rangle = \mathcal{T}$.
We say that $E\in\mathcal{T}$ is \textbf{finitely built} by $\mathcal{S}$ if $E\in\langle\mathcal{S}\rangle$.

Suppose that $\mathcal{T}$ is a compactly generated triangulated category (to clarify, we require that $\mathcal{T}$ admits all small coproducts in the definition of being compactly generated).
Denote the collection of compact objects of $\mathcal{T}$ by $\mathcal{T}^c$.
We say $G\in \mathcal{T}^c$ is a \textbf{compact generator} for $\mathcal{T}$ if for any $E\in\mathcal{T}$, one has $\operatorname{Hom}(G,E[n]) = 0$ for all $n\in \mathbb{Z}$ if, and only if, $E\cong 0$.
Classical generators for $\mathcal{T}^c$ coincide with compact generators for $\mathcal{T}$ \cite[\href{https://stacks.math.columbia.edu/tag/09SR}{Tag 09SR}]{StacksProject}.

\begin{example}
	\label{ex:perfect_compacts_classical_generator}
	Let $X$ be a quasi-compact quasi-separated scheme.
	Then the subcategory of compact objects in $D_{\operatorname{qc}}(X)$ coincides with $\operatorname{Perf}(X)$, and $D_{\operatorname{qc}}(X)$ is compactly generated by a single object of $\operatorname{Perf}(X)$.
	See \cite[Theorem 3.1.1]{Bondal/VandenBergh:2003}.
\end{example}

A triangulated subcategory of $\mathcal{T}$ is called \textbf{localizing} if it is closed under small coproducts.
Given a set $\mathcal{S}$ of objects in $\mathcal{T}$, we denote by $\overline{\langle \mathcal{S} \rangle}$ the smallest localizing subcategory of $\mathcal{T}$ which contains $\mathcal{S}$.

\begin{lemma}
	\label{lem:localizing_iff_cpt_gen}
	Assume $\mathcal{T}$ is compactly generated.
	Then $G\in \mathcal{T}^c$ compactly generates $\mathcal{T}$ if, and only if, $\overline{\langle G \rangle} = \mathcal{T}$.
\end{lemma}

\begin{proof}
	By \cite[Remark 1.16, Theorem 1.17, and Proposition 9.1.19]{Neeman:2001}, the Verdier localization $\pi \colon \mathcal{T}\to \mathcal{T}/\overline{\langle G \rangle}$ admits a right adjoint $Q\colon \mathcal{T}/\overline{\langle G \rangle}\to \mathcal{T}$.
	Hence, for every $E\in \mathcal{T}$, there exists a distinguished triangle
	\begin{displaymath}
		S_E \to E \to (Q\circ \pi)(E) \to S_E [1]
	\end{displaymath}
	with $S_E \in \overline{\langle G \rangle}$ and $(Q\circ \pi)(E)\in \overline{\langle G \rangle}^\perp$ (i.e.\ objects $A\in \mathcal{T}$ such that $\operatorname{Hom}(B,A)= 0$ for all $B\in \overline{\langle G \rangle}$).
	Therefore, the claim follows.
\end{proof}

\begin{lemma}
	\label{lem:equivalence_via_compacts}
	Suppose $F\colon \mathcal{T}\to \mathcal{S}$ is an exact functor between triangulated categories which are compactly generated by a single object.
	Assume $F$ has a right adjoint $G$ that preserves small coproducts.
	Then $F$ is an equivalence if, and only if, $F$ restricts to an equivalence $\mathcal{T}^c \to \mathcal{S}^c$.
\end{lemma}

\begin{proof}
	By \cite[Theorem 5.1]{Neeman:1996}, $F$ preserves compacts.
	Recall that an equivalence $\mathcal{T}\to\mathcal{S}$ induces an equivalence on compact objects $\mathcal{T}^c \to \mathcal{S}^c$.
	Conversely, if $F$ restricts to an equivalence $\mathcal{T}^c \to \mathcal{S}^c$, then this restriction has a right adjoint.
	Hence, \cite[Lemma 2.6]{Balmer/Dell'Ambrogio/Sanders:2016} implies $G$ preserves compact objects, and the adjunction of $F$ and $G$ on ambient categories restricts to an adjunction on categories of compact objects.
	To show $F$ is an equivalence, we show that $F$ and $G$ are both fully faithful.
	Denote the counit and unit of the adjunction respectively by $\epsilon$ and $\eta$.
	The strictly full subcategories $\mathcal{T}^\prime\subseteq\mathcal{T}$ and $\mathcal{S}^\prime\subseteq\mathcal{S}$ consisting respectively of objects where $\eta$ and $\epsilon$ are isomorphisms contain $\mathcal{T}^c$ and $\mathcal{S}^c$.
	One may verify that $\mathcal{T}^\prime$ and $\mathcal{S}^\prime$ are closed under shifts, iterated cones, and small coproducts (hence, homotopy colimits).
	By \cite[\href{https://stacks.math.columbia.edu/tag/09SN}{Tag 09SN}]{StacksProject}, each object in $\mathcal{T}$ and $\mathcal{S}$ can be obtained respectively from $\mathcal{T}^c$ and $\mathcal{S}^c$ using these operations.
	Hence, $\mathcal{T} \subseteq \mathcal{T}^\prime$ and $\mathcal{S} \subseteq \mathcal{S}^\prime$, which completes the proof.
\end{proof}

\subsection{$t$-structures}
\label{sec:prelim_t_structures}

We recall some notions for $t$-structures.
See \cite{Beilinson/Berstein/Deligne/Gabber:2018, Keller/Vossieck:1988}.
Let $\mathcal{T}$ be a triangulated category.
A pair of strictly full subcategories $\tau = (\mathcal{T}^{\leq 0}, \mathcal{T}^{\geq 0})$ of $\mathcal{T}$ is a \textbf{$t$-structure} if:
\begin{itemize}
	\item $\operatorname{Hom}(A,B) = 0$ for all $A\in\mathcal{T}^{\leq 0}$ and $B\in\mathcal{T}^{\geq 0}[-1]$,
	\item $\mathcal{T}^{\leq 0}[1]\subseteq\mathcal{T}^{\leq 0}$ and $\mathcal{T}^{\geq 0}[-1]\subseteq\mathcal{T}^{\geq 0}$,
	\item for every $E\in\mathcal{T}$, there is a distinguished triangle
	\begin{displaymath}
		\tau^{\leq 0}E \to E \to \tau^{\geq 1}E \to (\tau^{\leq 0}E)[1]
	\end{displaymath}
	with $\tau^{\leq 0}E\in\mathcal{T}^{\leq 0}$ and $\tau^{\geq 1}E\in\mathcal{T}^{\geq 0}[-1]$.
\end{itemize}
For $n\in\mathbb{Z}$, set $\mathcal{T}^{\leq n} := \mathcal{T}^{\leq 0}[-n]$ and $\mathcal{T}^{\geq n}:=\mathcal{T}^{\geq 0}[-n]$.
Denote by $\tau^{\leq n}$ and $\tau^{\geq n}$ the corresponding truncation functors.
The \textbf{$n$-th cohomology functor} with respect to $\tau$ is given by $H^n_\tau(E) := \tau^{\geq 0}\tau^{\leq 0}(E[n])$, and for each $k\in\mathbb{Z}$ there is a distinguished triangle
\begin{displaymath}
	\tau^{\leq k-1} E \to \tau^{\leq k} E \to H^k_\tau (E) [-k] \to (\tau^{\leq k-1} E)[1].
\end{displaymath}
See \cite[\S 1.3]{Beilinson/Berstein/Deligne/Gabber:2018}.
We say $\tau$ is \textbf{nondegenerate} if $\cap_{n\in\mathbb{Z}}\mathcal{T}^{\geq n} = \cap_{n\in\mathbb{Z}}\mathcal{T}^{\leq n} = \operatorname{add}(0)$ (i.e.\ both consist of zero objects), and \textbf{bounded} if for every $E\in\mathcal{T}$ there is an $n\geq 0$ such that $E[n]\in\mathcal{T}^{\leq 0}$ and $E[-n]\in\mathcal{T}^{\geq 0}$.

\begin{lemma}
	\label{lem:bounded_t_structure_is_nondegenerate}
	Let $\mathcal{T}$ be a triangulated category.
	A bounded $t$-structure $\tau = (\mathcal{T}^{\leq 0}, \mathcal{T}^{\geq 0})$ on $\mathcal{T}$ is nondegenerate.
\end{lemma}

\begin{proof}
	This is standard but we add it for convience. 
	Let $E\in \cap_{n\in \mathbb{Z}} \mathcal{T}^{\geq n}$.
	Assume $E$ is not the zero object.
	Since $\tau$ is bounded, there is an $N\geq 0$ such that $E[N]\in \mathcal{T}^{\leq 0}$ and $E[-N]\in \mathcal{T}^{\geq 0}$.
	This means $E\in \mathcal{T}^{\leq N}\cap \mathcal{T}^{\geq -N}$.
	However, $E\in \mathcal{T}^{\geq N+1} = \mathcal{T}^{\geq N}[-1]$, and so, $E\in \mathcal{T}^{\leq N}\cap \mathcal{T}^{\geq N+1}$.
	Since $\operatorname{Hom}(A,B) = 0$ for all $A \in \mathcal{T}^{\leq N}$ and $B \in \mathcal{T}^{\geq N}[-1]$, we have a contradiction.
	A similar argument shows $\cap_{n\in \mathbb{Z}} \mathcal{T}^{\leq n} = \operatorname{add}(0)$.
\end{proof}

\begin{lemma}
	\label{lem:bounded_t_structure_finitely_built_by_cohomology}
	Let $\mathcal{T}$ be a triangulated category equipped with a bounded $t$-structure $\tau$.
	Denote by $H^n_\tau$ the $n$-th cohomology functor with respect to $\tau$.
	Then any object $E$ in $\mathcal{T}$ is finitely built by $\oplus_{n\in \mathbb{Z}} H^n_\tau (E)$.
\end{lemma}

\begin{proof}
	This is standard but we add it for convience. 
	Fix $E\in \mathcal{T}$.
	We prove the claim by induction using the truncation triangles
	\begin{displaymath}
		\tau^{\leq k-1} E \to \tau^{\leq k} E \to H^k_\tau (E) [-k] \to (\tau^{\leq k-1} E)[1].
	\end{displaymath}
	By \Cref{lem:bounded_t_structure_is_nondegenerate}, $\tau$ is nondegenerate.
	This implies that the zero objects coincide with those objects whose $s$-th cohomology vanishes for all $s\in \mathbb{Z}$.
	Hence, $\langle \{ H^n_\tau (E) \}_{n\in \mathbb{Z}} \rangle_0$ consists of only zero objects.
	Thus, we can impose $E\not\cong 0$. 
	Set $k_1 = \min \{n\in \mathbb{Z} : H^n_\tau (E)\not\cong 0\}$. 
	Then $\tau^{\leq k_1 -1} E\cong 0$,
	and so, $\tau^{\leq k_1} E \cong H^{k_1}_\tau (E) [-k_1]$.
	Assume by induction that there exists an $N\geq0$ such that for all $0\leq s \leq N$ we have $\tau^{\leq k_1 +s} E \in \langle \{ H^n_\tau (E) \}_{n\in \mathbb{Z}} \rangle_{s+1}$.
	Consider the truncation triangle
	\begin{displaymath}
		\begin{aligned}
			\tau^{\leq k_1 + (N+1) -1} E & \to \tau^{\leq k_1 + (N+1) } E
			\\& \to H^{k_1 + (N+1)}_\tau (E) [- (k_1 + N+1)] \to (\tau^{\leq k_1 + (N+1) -1} E)[1].
		\end{aligned}
	\end{displaymath}
	By the inductive hypothesis,
	\begin{displaymath}
		\tau^{\leq k_1 + (N+1) -1} E\in \langle \{ H^n_\tau (E) \}_{n\in \mathbb{Z}} \rangle_{N+1},
	\end{displaymath}
	which implies that $\tau^{\leq k_1 + (N+1) } E\in \langle \{ H^n_\tau (E) \}_{n\in \mathbb{Z}} \rangle_{N+2}$.
	Since $\tau$ is a bounded $t$-structure, $H^n_\tau (E)\cong 0$ for $|n| \gg 0$, and hence $E \cong \tau^{\leq k_1 + L} E$ for $L\gg 0$.
	This completes the proof.
\end{proof}

\subsection{Schemes}
\label{sec:prelim_schemes}

We follow the conventions of \cite{StacksProject}.
Let $X$ be a scheme.
$D(X):=D(\operatorname{Mod}(X))$ is the derived category of $\mathcal{O}_X$-modules. 
$D_{\operatorname{qc}}(X)$ is the strictly full subcategory of $D(X)$ consisting of complexes with quasi-coherent cohomology.
$\operatorname{Perf}(X)$ is the strictly full subcategory of $D_{\operatorname{qc}}(X)$ consisting of the perfect complexes on $X$.
If $X$ is Noetherian, then $D_{\operatorname{coh}}^b(X)$ (resp.\ $D_{\operatorname{qc}}^b(X)$) is the strictly full subcategory of $D(X)$ consisting of complexes having bounded and coherent (resp.\ quasi-coherent) cohomology.
If $X$ is Noetherian, then $D^-_{\operatorname{coh}}(X)$ is the strictly full subcategory of $D_{\operatorname{qc}}(X)$ consisting of pseudocoherent complexes.
If $X$ is quasi-compact quasi-separated, then the standard $t$-structure restricts to a bounded $t$-structure on $D^b_{\operatorname{qc}}(X)$, whose $n$-th cohomology functor sends $E$ to $\mathcal{H}^n(E)$ (placed in degree zero).
If $X$ is affine, we might abuse notation and write $D^b_{\operatorname{coh}}(R):=D^b_{\operatorname{coh}}(X)$ where $R$ is the ring of global sections of $\mathcal{O}_X$.
Similar conventions will occur for the other categories.

Consider a morphism $f\colon Y \to X$ of Noetherian schemes.
The functor $\mathbf{R}f_\ast \colon D_{\operatorname{qc}}(Y)\to D_{\operatorname{qc}}(X)$ admits a right adjoint $f^\times$ \cite[\href{https://stacks.math.columbia.edu/tag/0A9E}{Tag 0A9E}]{StacksProject}.
If $f$ is separated and of finite type, then $f^!\colon D^+_{\operatorname{qc}}(X) \to D^+_{\operatorname{qc}}(Y)$ denotes the upper shriek functor \cite[\href{https://stacks.math.columbia.edu/tag/0A9Y}{Tag 0A9Y}]{StacksProject}.
If moreover $f$ is proper, then $f^\times$ and $f^!$ agree on $D^+_{\operatorname{qc}}(X)$ \cite[\href{https://stacks.math.columbia.edu/tag/0AU3}{Tag 0AU3}]{StacksProject}.
Recall that $\mathbf{L}f^\ast D^-_{\operatorname{coh}}(X)\subseteq D^-_{\operatorname{coh}}(Y)$ \cite[\href{https://stacks.math.columbia.edu/tag/08H4}{Tag 08H4}]{StacksProject}.
Denote by $\operatorname{\mathbf{R}\mathcal{H}\! \mathit{om}}$ the internal Hom functor for $D_{\operatorname{qc}}(X)$, and by $\operatorname{\mathbb{R}\mathcal{H}\! \mathit{om}}$ the internal Hom functor for $D(X)$.

Given $p\in X$, denote by $\sigma_p\colon \operatorname{Spec}(\mathcal{O}_{X,p})\to X$ and $i_p \colon \operatorname{Spec}(\kappa(p))\to X$ the natural morphisms.
The \textbf{support} of $E\in D_{\operatorname{qc}}(X)$ is the set $\operatorname{Supp}(E)$ of $p\in X$ such that $\mathbf{L}\sigma_p^\ast E\not\cong 0$.
Therefore,
\begin{displaymath}
	\operatorname{Supp}(E)
	=\bigcup_{i\in \mathbb{Z}}
	\operatorname{Supp}(\mathcal{H}^i(E)).
\end{displaymath}
If $X$ is Noetherian and $E\in D^-_{\operatorname{coh}}(X)$, then $\operatorname{Supp}(E)$ coincides with the set of $p\in X$ such that $\mathbf{L}i^\ast_p E\not\cong 0$ \cite[Lemma A.1]{Iyengar/Lipman/Neeman:2015}.
If furthermore $E\in D_{\operatorname{coh}}^b(X)$, then $\operatorname{Supp}(E)$ is closed \cite[\href{https://stacks.math.columbia.edu/tag/01BA}{Tag 01BA}]{StacksProject}.

\begin{lemma}
	\label{lem:faithfully_flat_conservative}
	Let $f\colon Y\to X$ be a faithfully flat morphism of Noetherian schemes.
	If $E\in D^b_{\operatorname{coh}}(X)$ satisfies $\mathbf{L}f^\ast E\cong 0$, then $E\cong 0$.
\end{lemma}

\begin{proof}
	Let $j\in\mathbb{Z}$.
	As $f$ is flat, $f^\ast \mathcal{H}^j(E) \cong \mathcal{H}^j(\mathbf{L}f^\ast E) \cong 0$.
	Since $f$ is flat and surjective, $\operatorname{Supp}(\mathcal{H}^j(E)) = f(\operatorname{Supp}(f^\ast\mathcal{H}^j(E)))$ \cite[\href{https://stacks.math.columbia.edu/tag/056J}{Tag 056J}]{StacksProject}, and so $\mathcal{H}^j(E) = 0$ for all $j\in\mathbb{Z}$.
\end{proof}

\subsection{Integral transforms}
\label{sec:prelim_integral_transforms}

Suppose $f_1 \colon Y_1 \to S$ and $f_2 \colon Y_2 \to S$ are morphisms of finite type to a Noetherian scheme.
Consider the fibered square
\begin{displaymath}
	\begin{tikzcd}
		{Y_1\times_S Y_2} & {Y_2} \\
		{Y_1} & {S.}
		\arrow["{f_1^\prime}", from=1-1, to=1-2]
		\arrow["{f_2^\prime}"', from=1-1, to=2-1]
		\arrow["{f_2}", from=1-2, to=2-2]
		\arrow["{f_1}"', from=2-1, to=2-2]
	\end{tikzcd}
\end{displaymath}
For any $K\in D_{\operatorname{qc}}(Y_1\times_S Y_2)$, the \textbf{integral transform with kernel} $K$ is the exact functor $\Phi_K\colon D_{\operatorname{qc}}(Y_1)\to D_{\operatorname{qc}}(Y_2)$ given by $E \mapsto \mathbf{R}(f_1^\prime)_\ast (\mathbf{L}(f_2^\prime)^\ast E \otimes^\mathbf{L} K)$.
We say $Y_1$ and $Y_2$ are \textbf{Fourier--Mukai $S$-partners} if there exists $K\in D^b_{\operatorname{coh}}(Y_1\times_S Y_2)$ such that $\Phi_K$ restricts to an equivalence $D^b_{\operatorname{coh}}(Y_1) \to D^b_{\operatorname{coh}}(Y_2)$.
When clear from context, the prefix `$S$-' is dropped.

Let $f\colon Y\to X$ be a proper morphism of Noetherian schemes.
An object $E\in D^-_{\operatorname{coh}}(Y)$ is called \textbf{$f$-perfect} if $\mathbf{R}f_\ast (E\otimes^{\mathbf{L}} \operatorname{Perf}(Y))\subseteq \operatorname{Perf}(X)$.
Now, suppose each $f_i$ is proper and flat.
An object $K\in D^-_{\operatorname{coh}}(Y_1\times_S Y_2)$ is called \textbf{relatively perfect over $Y_1$} (resp.\ \textbf{over $Y_2$}) if it is $f_2^\prime$-perfect (resp.\ $f_1^\prime$-perfect).
By \cite[Propositions 3.2 and 3.6]{Dutta/Lank/ManaliRahul:2025}, $K$ is relatively perfect over $Y_1$ if, and only if, $\Phi_K (D^b_{\operatorname{coh}}(Y_1))\subseteq D^b_{\operatorname{coh}}(Y_2)$, and $K$ is relatively perfect over $Y_2$ if, and only if, $\Phi_K (\operatorname{Perf}(Y_1))\subseteq \operatorname{Perf}(Y_2)$.
By symmetry, the analogous statements hold for integral transforms $D_{\operatorname{qc}}(Y_2)\to D_{\operatorname{qc}}(Y_1)$ with kernels in $D^-_{\operatorname{coh}}(Y_1\times_S Y_2)$.

\begin{lemma}
	\label{lem:relatively_perfect_pseudocoherent}
	Let $f\colon Y \to X$ be a proper morphism of Noetherian schemes.
	Then the following are equivalent for any $E\in D^-_{\operatorname{coh}} (Y)$:
	\begin{enumerate}
		\item \label{lem:relatively_perfect_pseudocoherent1} $E$ is $f$-perfect
		\item \label{lem:relatively_perfect_pseudocoherent2} $E \otimes^{\mathbf{L}} \mathbf{L}f^\ast D^b_{\operatorname{coh}}(X)\subseteq D^b_{\operatorname{coh}}(Y)$
		\item \label{lem:relatively_perfect_pseudocoherent3} $E \otimes^{\mathbf{L}} \mathbf{L}f^\ast D^b_{\operatorname{qc}}(X)\subseteq D^b_{\operatorname{qc}}(Y)$.
	\end{enumerate}
	In particular, any $f$-perfect object of $D^-_{\operatorname{coh}}(Y)$ belongs to $D^b_{\operatorname{coh}}(Y)$.
\end{lemma}

\begin{proof}
	The equivalence $\eqref{lem:relatively_perfect_pseudocoherent1} \iff \eqref{lem:relatively_perfect_pseudocoherent2}$ can be argued affine locally on the target by using \cite[Lemma 3.4]{Dutta/Lank/ManaliRahul:2025}, whereas $\eqref{lem:relatively_perfect_pseudocoherent1} \iff \eqref{lem:relatively_perfect_pseudocoherent3}$ is \cite[Theorem 4.3]{AlonsoTarrio/JeremiasLopez/SanchodeSalas:2023}.
	The last claim follows from \eqref{lem:relatively_perfect_pseudocoherent2} applied to $\mathcal{O}_X$.
\end{proof}

We work in the following setup.

\begin{setup}
	\label{setup:fm_cube}
	Let $t\colon T \to S$ be a morphism of Noetherian schemes.
	Suppose $f_1\colon Y_1 \to S$ and $f_2\colon Y_2 \to S$ are proper flat morphisms.
	Consider the commutative diagram
	\begin{equation}
		\label{diag:fm_cube}
		\begin{tikzcd}
			{Y_1\times_S Y_2\times_S T} && {Y_2\times_S T} \\
			& {Y_1\times_S T} &&& T \\
			{Y_1\times_S Y_2} && {Y_2} \\
			& {Y_1} &&& S
			\arrow["{g_1^\prime}"', from=1-1, to=1-3]
			\arrow["{g_2^\prime}"', from=1-1, to=2-2]
			\arrow["{t^\prime}"', from=1-1, to=3-1]
			\arrow["{g_2}", from=1-3, to=2-5]
			\arrow["{t_2}"'{pos=0.8}, dashed, from=1-3, to=3-3]
			\arrow["{g_1}", from=2-2, to=2-5]
			\arrow["{t_1}"{pos=0.7}, from=2-2, to=4-2]
			\arrow["t", from=2-5, to=4-5]
			\arrow["{f^\prime_1}"{pos=0.3}, dashed, from=3-1, to=3-3]
			\arrow["{f^\prime_2}"', from=3-1, to=4-2]
			\arrow["{f_2}", from=3-3, to=4-5]
			\arrow["{f_1}"', from=4-2, to=4-5]
		\end{tikzcd}
	\end{equation}
	obtained by base change along $t$.
\end{setup}

\begin{lemma}
	\label{lem:fm_cube_tor_independent}
	Consider \Cref{setup:fm_cube}.
	Then each face of \eqref{diag:fm_cube} is a tor-independent square in the sense of \cite[Definition 3.10.2]{Lipman:2009}.
\end{lemma}

\begin{proof}
	As each $f_i$, $f^\prime_i$, $g_i$, and $g^\prime_i$ is flat, \cite[Remark 22.94]{Gortz/Wedhorn:2023} yields the desired claim.
\end{proof}

\begin{remark}
	\label{rem:fm_cube_separated}
	Consider \Cref{setup:fm_cube}.
	If $S$ is separated and $t$ is affine, then every scheme appearing in \eqref{diag:fm_cube} is separated.
	Indeed, proper morphisms and affine morphisms are separated, and separatedness is stable under base change and composition.
\end{remark}

\begin{lemma}
	\label{lem:fm_cube_flat_base_change}
	Consider \Cref{setup:fm_cube}.
	Let $K\in D_{\operatorname{qc}}(Y_1\times_S Y_2)$.
	Then on $D_{\operatorname{qc}}$ there is a natural isomorphism
	\begin{displaymath}
		\alpha^K \colon \mathbf{L} t_2^\ast \circ \Phi_{K} \to \Phi_{\mathbf{L}(t^\prime)^\ast K} \circ \mathbf{L} t_1^\ast.
	\end{displaymath}
\end{lemma}

\begin{proof}
	For each $E\in D_{\operatorname{qc}}(Y_1)$, there is a string of natural isomorphisms:
	\begin{displaymath}
		\begin{aligned}
			\mathbf{L} t_2^\ast \Phi_{K} (E)
			&= \mathbf{L} t_2^\ast \mathbf{R}(f^\prime_1)_\ast (\mathbf{L} (f^\prime_2)^\ast E \otimes^{\mathbf{L}} K)
			\\&\cong \mathbf{R}(g^\prime_1)_\ast \mathbf{L}(t^\prime)^\ast (\mathbf{L} (f^\prime_2)^\ast E \otimes^{\mathbf{L}} K)
			&& \textrm{(flat base change)}
			\\&\cong \mathbf{R}(g^\prime_1)_\ast ( \mathbf{L}(t^\prime)^\ast  \mathbf{L} (f^\prime_2)^\ast E \otimes^{\mathbf{L}}  \mathbf{L}(t^\prime)^\ast K)
			&& \textrm{(monoidality)}
			\\&\cong \mathbf{R}(g^\prime_1)_\ast ( \mathbf{L} (g^\prime_2)^\ast  \mathbf{L} t_1^\ast E \otimes^{\mathbf{L}}  \mathbf{L}(t^\prime)^\ast K)
			&& \textrm{(pseudofunctoriality)}
			\\&= \Phi_{\mathbf{L}(t^\prime)^\ast K}( \mathbf{L} t_1^\ast E).
		\end{aligned}
	\end{displaymath}
	This completes the proof.
\end{proof}

\section{Relative perfectness under base change}
\label{sec:relative_perfectness}

\subsection{Relative perfectness and finite tor-dimension}
\label{sec:relative_perfectness_tor_dimension}

\begin{lemma}
	\label{lem:ALS_implies_finite_flat_dim}
	Let $f\colon Y \to X$ be a morphism of quasi-compact quasi-separated schemes.
	For any $E\in D^b_{\operatorname{qc}}(Y)$, the following are equivalent:
	\begin{enumerate}
		\item \label{lem:ALS_implies_finite_flat_dim1} $E\otimes^{\mathbf{L}} \mathbf{L}f^\ast  D^b_{\operatorname{qc}}(X) \subseteq D^b_{\operatorname{qc}}(Y)$
		\item \label{lem:ALS_implies_finite_flat_dim2} there exists $[a,b]\subseteq \mathbb{Z}$ such that $\mathcal{H}^j (E\otimes^{\mathbf{L}} \mathbf{L} f^\ast M)\cong 0$ for all $j\not\in [a,b]$ and $M\in \operatorname{Qcoh}(X)$.
	\end{enumerate}
\end{lemma}

\begin{proof}
	First, we show $\eqref{lem:ALS_implies_finite_flat_dim1} \implies \eqref{lem:ALS_implies_finite_flat_dim2}$.
	Assume the contrary.
	Then for each $n \geq 1$ there exists an $M_n\in \operatorname{Qcoh}(X)$ such that $\mathcal{H}^j (E\otimes^{\mathbf{L}} \mathbf{L} f^\ast M_n)\not\cong 0$ for some $j\not\in [-n,n]$.
	Set $M:= \oplus_{n\geq 1} M_n$.
	By hypothesis, $E\otimes^{\mathbf{L}} \mathbf{L}f^\ast  M \in D^b_{\operatorname{qc}}(Y)$, and hence there exists $[a,b]\subseteq \mathbb{Z}$ such that $\mathcal{H}^j (E\otimes^{\mathbf{L}} \mathbf{L} f^\ast M)\cong 0$ for all $j\not\in [a,b]$.
	Since $\mathcal{H}^j (E\otimes^{\mathbf{L}} \mathbf{L} f^\ast M_n)$ is a direct summand of $\mathcal{H}^j (E\otimes^{\mathbf{L}} \mathbf{L} f^\ast M)$ for each $n\geq 1$, $\mathcal{H}^j (E\otimes^{\mathbf{L}} \mathbf{L} f^\ast M_n)\cong 0$ for all $j\not\in [a,b]$ and $n\geq 1$.
	This leads to a contradiction.
	Indeed, for $n> \max\{|a|,|b|\}$ there exists $j\not\in [-n,n]$ with
	$\mathcal{H}^j(E\otimes^{\mathbf{L}} \mathbf{L}f^\ast M_n)\not\cong 0$.

	Next, we prove that $\eqref{lem:ALS_implies_finite_flat_dim2} \implies \eqref{lem:ALS_implies_finite_flat_dim1}$. 
	Let $B\in D^b_{\operatorname{qc}}(X)$.
	By \Cref{lem:bounded_t_structure_finitely_built_by_cohomology} applied to the standard $t$-structure on $D^b_{\operatorname{qc}}(X)$, $B$ is finitely built by its cohomology sheaves $\mathcal{H}^i (B)$.
	Moreover, the hypothesis implies $E\otimes^{\mathbf{L}} \mathbf{L}f^\ast \mathcal{H}^i (B) \in D^b_{\operatorname{qc}}(Y)$ for all $i\in \mathbb{Z}$.
	Since $B\in \langle \{ \mathcal{H}^i (B) \}_{i\in\mathbb{Z}} \rangle$ and $D^b_{\operatorname{qc}}(Y)$ is a thick subcategory of $D_{\operatorname{qc}}(Y)$, it follows that
	\begin{displaymath}
		E\otimes^{\mathbf{L}} \mathbf{L} f^\ast B\in \langle \{ E\otimes^{\mathbf{L}} \mathbf{L} f^\ast \mathcal{H}^i (B) \}_{i\in \mathbb{Z}} \rangle\subseteq D^b_{\operatorname{qc}}(Y).
	\end{displaymath}
	This completes the proof.
\end{proof}

\begin{proposition}
	\label{prop:relative_perf_for_schemes}
	Let $f\colon Y \to X$ be a morphism of quasi-compact quasi-separated schemes.
	For any $E\in D^b_{\operatorname{qc}}(Y)$, the following are equivalent:
	\begin{enumerate}
		\item \label{prop:relative_perf_for_schemes1} $E\otimes^{\mathbf{L}} \mathbf{L}f^\ast  D^b_{\operatorname{qc}}(X) \subseteq D^b_{\operatorname{qc}}(Y)$
		\item \label{prop:relative_perf_for_schemes2} $E$ has finite tor-dimension as an object of $D(f^{-1}\mathcal{O}_X)$.
	\end{enumerate}
\end{proposition}

\begin{proof}
	First, we show $\eqref{prop:relative_perf_for_schemes1} \implies \eqref{prop:relative_perf_for_schemes2}$.
	By \Cref{lem:ALS_implies_finite_flat_dim}, there exists $[a,b]\subseteq \mathbb{Z}$ such that $\mathcal{H}^j (E\otimes^{\mathbf{L}} \mathbf{L} f^\ast M)\cong 0$ for all $j\not\in [a,b]$ and $M\in \operatorname{Qcoh}(X)$.
	Following \cite[pg.\ 78, Example 2.6.7]{Lipman:2009}, this condition is equivalent to $E$ having finite flat $f$-amplitude in $[a,b]$ (see loc.\ cit.\ for the definition).
	Moreover, this condition is equivalent to the requirement that for each $p\in Y$, $E_p$ is isomorphic in $D(\mathcal{O}_{X,f(p)})$ to a complex of flat $\mathcal{O}_{X,f(p)}$-modules vanishing in degrees outside of $[a,b]$; see the comments above \cite[Eq.\ 2.7.6.1]{Lipman:2009}.
	Applying \cite[Proposition 3.3]{Illusie:1971}, this stalk local condition is equivalent to $E$ having finite tor-dimension in $[a,b]$ as an object of $D(f^{-1}\mathcal{O}_X)$ (see also \cite[Proposition 21.169 and Lemma 21.171]{Gortz/Wedhorn:2023}).
	It follows that $\eqref{prop:relative_perf_for_schemes1} \implies \eqref{prop:relative_perf_for_schemes2}$.
	Conversely, the same references show that \eqref{prop:relative_perf_for_schemes2} implies \Cref{lem:ALS_implies_finite_flat_dim}\eqref{lem:ALS_implies_finite_flat_dim2}.
	Thus, $\eqref{prop:relative_perf_for_schemes2}\implies \eqref{prop:relative_perf_for_schemes1}$.
\end{proof}

\begin{corollary}
	\label{cor:relative_perfectness_tor_dimension}
	Let $f\colon Y\to X$ be a proper morphism of Noetherian schemes.
	Then $E\in D^-_{\operatorname{coh}}(Y)$ is $f$-perfect if, and only if, $E$ has finite tor-dimension as an object of $D(f^{-1}\mathcal{O}_X)$.
	In particular, if $f$ is moreover flat, then the $f$-perfect objects of $D^-_{\operatorname{coh}}(Y)$ are precisely the $X$-perfect objects in the sense of \cite[\href{https://stacks.math.columbia.edu/tag/0DI0}{Tag 0DI0}]{StacksProject}.
\end{corollary}

\begin{proof}
	If $E$ is $f$-perfect, then \Cref{lem:relatively_perfect_pseudocoherent} implies $E\in D^b_{\operatorname{coh}}(Y)$ and $E\otimes^{\mathbf{L}}\mathbf{L}f^\ast D^b_{\operatorname{qc}}(X)\subseteq D^b_{\operatorname{qc}}(Y)$.
	Hence, by \Cref{prop:relative_perf_for_schemes}, $E$ has finite tor-dimension over $f^{-1}\mathcal{O}_X$.
	Conversely, if $E$ has tor-amplitude in $[a,b]$ as an object of $D(f^{-1}\mathcal{O}_X)$, then $\mathcal{H}^i(E)\cong\mathcal{H}^i(E\otimes^{\mathbf{L}}_{f^{-1}\mathcal{O}_X} f^{-1}\mathcal{O}_X) = 0$ for $i\not\in[a,b]$, and so $E\in D^b_{\operatorname{qc}}(Y)$.
	Hence, \Cref{prop:relative_perf_for_schemes} and \Cref{lem:relatively_perfect_pseudocoherent} imply $E$ is $f$-perfect.
	For the last claim, note that $Y$ is quasi-compact, so $E$ has finite tor-dimension over $f^{-1}\mathcal{O}_X$ if, and only if, it does so locally on $Y$.
\end{proof}

\subsection{Base change}
\label{sec:relative_perfectness_base_change}

\begin{lemma}
	\label{lem:relative_perfectness_affine_base_change}
	Consider \Cref{setup:fm_cube} with $t$ an affine morphism.
	If $K\in D^b_{\operatorname{coh}}(Y_1\times_S Y_2)$ is relatively perfect over $Y_1$ (resp.\ over $Y_2$), then $\mathbf{L} (t^\prime)^\ast K$ is relatively perfect over $Y_1\times_S T$ (resp.\ over $Y_2\times_S T$).
	In particular, $\Phi_{\mathbf{L} (t^\prime)^\ast K}$ restricts to an exact functor on $D^b_{\operatorname{coh}}$ (resp.\ on $\operatorname{Perf}$).
\end{lemma}

\begin{remark}
	\Cref{lem:relative_perfectness_affine_base_change} can be seen directly from \cite[Proposition 3.11]{AlonsoTarrio/JeremiasLopez/SanchodeSalas:2023} together with \Cref{lem:relatively_perfect_pseudocoherent}, as the faces of \eqref{diag:fm_cube} are tor-independent (\Cref{lem:fm_cube_tor_independent}); in fact, loc.\ cit.\ does not require $t$ to be affine.
	We give a different proof for the sake of interest.
\end{remark}

\begin{proof}
	[Proof of \Cref{lem:relative_perfectness_affine_base_change}]
	We prove the case of $K$ being relatively perfect over $Y_2$; a similar argument applies over $Y_1$.
	By definition,
	\begin{displaymath}
		\mathbf{R}(f_1^\prime)_\ast (K \otimes^{\mathbf{L} }\operatorname{Perf}(Y_1\times_S Y_2)) \subseteq \operatorname{Perf}(Y_2).
	\end{displaymath}
	By base change, $t_1$, $t_2$, and $t^\prime$ are affine morphisms.
	Choose $G\in \operatorname{Perf}(Y_1\times_S Y_2)$ such that $\operatorname{Perf}(Y_1\times_S Y_2) = \langle G \rangle$.
	Then $\mathbf{L}(t^\prime)^\ast G$ satisfies
	\begin{displaymath}
		\operatorname{Perf}(Y_1\times_S Y_2\times_S T) = \langle \mathbf{L} (t^\prime)^\ast G \rangle,
	\end{displaymath}
	see \cite[\href{https://stacks.math.columbia.edu/tag/0BQT}{Tag 0BQT}]{StacksProject}.
	Also, by base change, $g_1$, $g_2$, $f^\prime_2$, $f^\prime_1$, $g_2^\prime$, and $g_1^\prime$ are proper flat morphisms.
	Using flat base change \cite[Remark 22.94 and Theorem 22.99]{Gortz/Wedhorn:2023}, it follows that
	\begin{displaymath}
		\begin{aligned}
			\mathbf{R} (g^\prime_1)_\ast & (\mathbf{L}(t^\prime)^\ast K \otimes^\mathbf{L} \operatorname{Perf}(Y_1\times_S Y_2 \times_S T))
			\\&\subseteq \mathbf{R} (g^\prime_1)_\ast \langle \mathbf{L}(t^\prime)^\ast (K \otimes^\mathbf{L} G) \rangle
			\\&\subseteq \langle \mathbf{R} (g^\prime_1)_\ast \mathbf{L}(t^\prime)^\ast (K \otimes^\mathbf{L} G) \rangle
			\\&\subseteq \langle \mathbf{L}t_2^\ast \mathbf{R}(f_1^\prime)_\ast (K \otimes^\mathbf{L} G) \rangle
			\\&\subseteq \operatorname{Perf}(Y_2\times_S T).
		\end{aligned}
	\end{displaymath}
	Moreover, $\mathbf{L}(t^\prime)^\ast K$ is pseudocoherent \cite[\href{https://stacks.math.columbia.edu/tag/08H4}{Tag 08H4}]{StacksProject}.
	Therefore, $\mathbf{L}(t^\prime)^\ast K$ is relatively perfect over $Y_2\times_S T$.
	The last claim follows from \cite[Propositions 3.2 and 3.6]{Dutta/Lank/ManaliRahul:2025}.
\end{proof}

\begin{lemma}
	\label{lem:pseudocoherence_faithfully_flat}
	Let $f\colon Y \to X$ be a faithfully flat morphism of Noetherian schemes.
	An object $E\in D_{\operatorname{qc}}(X)$ belongs to $D^-_{\operatorname{coh}}(X)$ (resp.\ $\operatorname{Perf}(X)$) if, and only if, $\mathbf{L}f^\ast E\in D^-_{\operatorname{coh}}(Y)$ (resp.\ $\mathbf{L}f^\ast E\in \operatorname{Perf}(Y)$).
\end{lemma}

\begin{proof}
	See \cite[Proposition 22.52]{Gortz/Wedhorn:2023}.
\end{proof}

\begin{lemma}
	\label{lem:perfect_around_a_point}
	Let $X$ be a scheme, $E\in D_{\operatorname{qc}}^b(X)$ be pseudocoherent, and $p\in X$.
	Then the following are equivalent:
	\begin{enumerate}
		\item \label{lem:perfect_around_a_point1} there is an open immersion $j\colon U \to X$ such that $p\in U$ and $\mathbf{L}j^\ast E$ is perfect
		\item \label{lem:perfect_around_a_point2} $\mathbf{L}\sigma_p^\ast E$ is perfect
		\item \label{lem:perfect_around_a_point3} $\mathbf{L}i_p^\ast E$ is bounded.
	\end{enumerate}
\end{lemma}

\begin{proof}
	To see that $\eqref{lem:perfect_around_a_point1}\implies \eqref{lem:perfect_around_a_point2}$ and $\eqref{lem:perfect_around_a_point2}\implies \eqref{lem:perfect_around_a_point3}$, use the natural factorizations
	\begin{displaymath}
		\operatorname{Spec}(\kappa(p))
		\to\operatorname{Spec}(\mathcal{O}_{X,p})
		\to U\xrightarrow{j} X
	\end{displaymath}
	of $\sigma_p$ and of $i_p$, and the fact that derived pullback preserves perfect complexes \cite[\href{https://stacks.math.columbia.edu/tag/09UA}{Tag 09UA}]{StacksProject}.
	Lastly, we check $\eqref{lem:perfect_around_a_point3}\implies \eqref{lem:perfect_around_a_point1}$.
	In this case, we may assume $X$ is affine, and the result follows from \cite[Theorem 2.3, (iii)$\implies$(i)]{AlonsoTarrio/JeremiasLopez/SanchodeSalas:2023}.
\end{proof}

\begin{proposition}
	\label{prop:perfection_is_stalk_local_property}
	Let $X$ be a quasi-compact quasi-separated scheme and $E\in D^b_{\operatorname{qc}}(X)$ be pseudocoherent.
	Then the following are equivalent:
	\begin{enumerate}
		\item \label{prop:perfection_is_stalk_local_property1} $E$ is perfect
		\item \label{prop:perfection_is_stalk_local_property2} $\mathbf{L}i_p^\ast E$ is bounded for every $p \in X$
		\item \label{prop:perfection_is_stalk_local_property3} $\mathbf{L}i_p^\ast E$ is bounded for every closed point $p \in X$
		\item \label{prop:perfection_is_stalk_local_property4} $\mathbf{L}\sigma_p^\ast E$ is perfect for every $p \in X$
		\item \label{prop:perfection_is_stalk_local_property5} $\mathbf{L}\sigma_p^\ast E$ is perfect for every closed point $p \in X$.
	\end{enumerate}
\end{proposition}

\begin{proof}
	It is straightforward to see that $\eqref{prop:perfection_is_stalk_local_property2}\implies \eqref{prop:perfection_is_stalk_local_property3}$.
	By \cite[Theorem 2.3]{AlonsoTarrio/JeremiasLopez/SanchodeSalas:2023}, we have $\eqref{prop:perfection_is_stalk_local_property1}\iff \eqref{prop:perfection_is_stalk_local_property2}$.
	Moreover, from \Cref{lem:perfect_around_a_point}, we have $\eqref{prop:perfection_is_stalk_local_property4} \iff \eqref{prop:perfection_is_stalk_local_property2}$ and $\eqref{prop:perfection_is_stalk_local_property5} \iff \eqref{prop:perfection_is_stalk_local_property3}$.
	Let $p\in X$.
	Since $X$ is quasi-compact, there is a closed point $x\in X$ such that $x\in\overline{\{p\}}$ \cite[Exercise 3.13]{Gortz/Wedhorn:2020}.
	Hence, $\sigma_p$ factors as the natural morphism $\operatorname{Spec}(\mathcal{O}_{X,p})\to\operatorname{Spec}(\mathcal{O}_{X,x})$ followed by $\sigma_x$, and so $\eqref{prop:perfection_is_stalk_local_property5} \iff \eqref{prop:perfection_is_stalk_local_property4}$ by \cite[\href{https://stacks.math.columbia.edu/tag/09UA}{Tag 09UA}]{StacksProject}.
\end{proof}

\begin{lemma}
	\label{lem:reflecting_relative_perfectness}
	Consider \Cref{setup:fm_cube}.
	Let $K\in D^b_{\operatorname{coh}}(Y_1\times_S Y_2)$.
	\begin{enumerate}
		\item \label{lem:reflecting_relative_perfectness1} If $t$ is affine and faithfully flat and $\mathbf{L}(t^\prime)^\ast K$ is relatively perfect over $Y_1\times_S T$ (resp.\ $Y_2\times_S T$), then $K$ is relatively perfect over $Y_1$ (resp.\ $Y_2$).
		\item \label{lem:reflecting_relative_perfectness2} If $\mathbf{L}(t^\prime)^\ast K$ is relatively perfect over $Y_1\times_S \operatorname{Spec}(k)$ (resp.\ $Y_2\times_S \operatorname{Spec}(k)$) for all morphisms $t\colon \operatorname{Spec}(k)\to S$ from a field, then $K$ is relatively perfect over $Y_1$ (resp.\ $Y_2$).
		\item \label{lem:reflecting_relative_perfectness3} If $\mathbf{L}(t^\prime)^\ast K$ is relatively perfect over $Y_1\times_S \operatorname{Spec}(\kappa(p))$ (resp.\ $Y_2\times_S \operatorname{Spec}(\kappa(p))$) for all closed immersions $t\colon \operatorname{Spec}(\kappa(p))\to S$ associated to a closed point $p\in S$, then $K$ is relatively perfect over $Y_1$ (resp.\ $Y_2$).
	\end{enumerate}
\end{lemma}

\begin{proof}
	We prove the case for relative perfectness over $Y_1$.
	To show the other case, one can argue as below and use \Cref{lem:pseudocoherence_faithfully_flat}.

	We prove the first claim.
	Let $E\in D^b_{\operatorname{coh}}(Y_1)$.
	As $t_1$ is faithfully flat, it follows that $\mathbf{L} t_1^\ast E\in D^b_{\operatorname{coh}}(Y_1\times_S T)$.
	By \Cref{lem:fm_cube_flat_base_change}, $\mathbf{L} t_2^\ast \Phi_K (E)\cong \Phi_{\mathbf{L}(t^\prime)^\ast K} (\mathbf{L} t^\ast_1 E)$.
	Our hypothesis ensures that $\Phi_{\mathbf{L}(t^\prime)^\ast K} (\mathbf{L} t_1^\ast E) \in D^b_{\operatorname{coh}}(Y_2\times_S T)$.
	From \Cref{lem:pseudocoherence_faithfully_flat}, it follows that $\Phi_K (E) \in D_{\operatorname{coh}}^-(Y_2)$.
	Hence, it suffices to show that $\mathcal{H}^j (\Phi_K (E))=0$ for all but finitely many $j\in \mathbb{Z}$.
	As $t_2$ is flat and surjective, for each $j\in \mathbb{Z}$,
	\begin{displaymath}
		\begin{aligned}
			\operatorname{Supp}(\mathcal{H}^j (\Phi_K (E)) )
			&= t_2 (\operatorname{Supp}(t_2^\ast \mathcal{H}^j (\Phi_K (E)) )) && \textrm{(\cite[\href{https://stacks.math.columbia.edu/tag/056J}{Tag 056J}]{StacksProject})}
			\\&= t_2 (\operatorname{Supp}(\mathcal{H}^j (\mathbf{L}t_2^\ast \Phi_K (E)) ))  && \textrm{(flatness of $t_2$)}.
		\end{aligned}
	\end{displaymath}
	Hence, $\Phi_K (E)$ has bounded cohomology.
	Consequently, $\Phi_K (D^b_{\operatorname{coh}}(Y_1)) \subseteq D^b_{\operatorname{coh}}(Y_2)$, which implies $K$ is relatively perfect over $Y_1$ by \cite[Propositions 3.2 and 3.6]{Dutta/Lank/ManaliRahul:2025}.

	Next, we check the second claim.
	Let $P\in \operatorname{Perf}(Y_1\times_S Y_2)$.
	Note that $\mathbf{R}(f^\prime_2)_\ast ( K \otimes^{\mathbf{L}} P)$ is bounded and pseudocoherent because $K \otimes^{\mathbf{L}} P$ is a bounded pseudocoherent complex and $f^\prime_2$ is proper.
	By \cite[Theorem 2.3]{AlonsoTarrio/JeremiasLopez/SanchodeSalas:2023}, it suffices to show that $\mathbf{L}i_p^\ast \mathbf{R}(f^\prime_2)_\ast ( K \otimes^{\mathbf{L}} P)$ is bounded for every $p\in Y_1$.
	Choose $p\in Y_1$.
	Denote by $t\colon \operatorname{Spec}(\kappa(f_1(p)))\to S$ the natural morphism.
	There is a commutative diagram
	\begin{displaymath}
		\begin{tikzcd}
			{\operatorname{Spec}(\kappa(p))} && \\
			& {Y_1\times_S \operatorname{Spec}(\kappa(f_1(p)))} & {\operatorname{Spec}(\kappa(f_1(p)))} \\
			& {Y_1} & S
			\arrow["h"{description}, from=1-1, to=2-2]
			\arrow["{h_2}", bend left=30, from=1-1, to=2-3]
			\arrow["{i_p}"', bend right=30, from=1-1, to=3-2]
			\arrow["{g_1}", from=2-2, to=2-3]
			\arrow["{t_1}"', from=2-2, to=3-2]
			\arrow["t", from=2-3, to=3-3]
			\arrow["{f_1}", from=3-2, to=3-3]
		\end{tikzcd}
	\end{displaymath}
	where $h$ is induced by the universal property of the fiber product.
	Using the hypothesis, it follows that
	\begin{displaymath}
		\begin{aligned}
			\mathbf{R}(g^\prime_2)_\ast
			& (\mathbf{L}(t^\prime)^\ast ( K \otimes^{\mathbf{L}} P ) )
			\\&\in \mathbf{R}(g^\prime_2)_\ast (\mathbf{L}(t^\prime)^\ast K \otimes^{\mathbf{L}} \operatorname{Perf}(Y_1\times_S Y_2 \times_S \operatorname{Spec}(\kappa(f_1(p)))))
			\\&\subseteq \operatorname{Perf}(Y_1\times_S \operatorname{Spec}(\kappa(f_1(p)))).
		\end{aligned}
	\end{displaymath}
	By tor-independent base change for the fibered square
	\begin{displaymath}
		\begin{tikzcd}
			{Y_1\times_S Y_2\times_S \operatorname{Spec}(\kappa(f_1(p)))} & {Y_1\times_S \operatorname{Spec}(\kappa(f_1(p)))} \\
			{Y_1\times_S Y_2} & {Y_1}
			\arrow["{g_2^\prime}", from=1-1, to=1-2]
			\arrow["{t^\prime}"', from=1-1, to=2-1]
			\arrow["{t_1}", from=1-2, to=2-2]
			\arrow["{f^\prime_2}"', from=2-1, to=2-2]
		\end{tikzcd}
	\end{displaymath}
	\cite[Remark 22.94 and Theorem 22.99]{Gortz/Wedhorn:2023}, we have
	\begin{displaymath}
		\mathbf{L}t_1^\ast \mathbf{R}(f^\prime_2)_\ast ( K \otimes^{\mathbf{L}} P)
		\cong \mathbf{R}(g^\prime_2)_\ast \mathbf{L}(t^\prime)^\ast ( K \otimes^{\mathbf{L}} P)
		\in \operatorname{Perf}(Y_1\times_S \operatorname{Spec}(\kappa(f_1(p)))).
	\end{displaymath}
	This implies that
	\begin{displaymath}
		\mathbf{L}i_p^\ast \mathbf{R}(f^\prime_2)_\ast ( K \otimes^{\mathbf{L}} P)
		\cong \mathbf{L}h^\ast \mathbf{L}t_1^\ast \mathbf{R}(f^\prime_2)_\ast ( K \otimes^{\mathbf{L}} P)
		\in D^b_{\operatorname{coh}}(\kappa(p)),
	\end{displaymath}
	because derived pullback preserves perfect complexes \cite[\href{https://stacks.math.columbia.edu/tag/09UA}{Tag 09UA}]{StacksProject}.
	Thus, the complex $\mathbf{R}(f^\prime_2)_\ast ( K \otimes^{\mathbf{L}} P)$ is perfect, and so $K$ is relatively perfect over $Y_1$.

	Lastly, we check the third claim.
	The argument above shows that $\mathbf{L}i_p^\ast \mathbf{R}(f^\prime_2)_\ast ( K \otimes^{\mathbf{L}} P)$ is bounded for every closed point $p\in Y_1$.
	Indeed, $f_1$ is proper, and hence closed on the underlying topological spaces, so $f_1(p)$ is a closed point of $S$.
	Consequently, the desired claim follows from \Cref{prop:perfection_is_stalk_local_property}.
\end{proof}

\begin{theorem}
	\label{thm:relative_perfectness_base_change}
	Consider \Cref{setup:fm_cube}.
	Let $K\in D^b_{\operatorname{coh}}(Y_1\times_S Y_2)$.
	Then the following are equivalent:
	\begin{enumerate}
		\item \label{thm:relative_perfectness_base_change1} $K$ is relatively perfect over $Y_1$ (resp.\ $Y_2$)
		\item \label{thm:relative_perfectness_base_change2} $\mathbf{L}(t^\prime)^\ast K$ is relatively perfect over $Y_1\times_S T$ (resp.\ $Y_2\times_S T$) for every affine morphism $t\colon T\to S$ from a Noetherian scheme
		\item \label{thm:relative_perfectness_base_change3} $\mathbf{L}(t^\prime)^\ast K$ is relatively perfect over $Y_1\times_S \operatorname{Spec}(k)$ (resp.\ $Y_2\times_S \operatorname{Spec}(k)$) for every morphism $t\colon \operatorname{Spec}(k)\to S$ from a field
		\item \label{thm:relative_perfectness_base_change4} $\mathbf{L}(t^\prime)^\ast K$ is relatively perfect over $Y_1\times_S \operatorname{Spec}(\kappa(p))$ (resp.\ $Y_2\times_S \operatorname{Spec}(\kappa(p))$) for every closed immersion $t\colon \operatorname{Spec}(\kappa(p))\to S$ associated to a closed point $p\in S$
		\item \label{thm:relative_perfectness_base_change5} $\mathbf{L}(t^\prime)^\ast K$ is relatively perfect over $Y_1\times_S T$ (resp.\ $Y_2\times_S T$) for some affine surjective morphism $t\colon T\to S$ from a Noetherian scheme.
	\end{enumerate}
\end{theorem}

\begin{proof}
	We prove the case of relative perfectness over $Y_1$; the other can be shown analogously.
	By \Cref{lem:relative_perfectness_affine_base_change}, \eqref{thm:relative_perfectness_base_change1} implies \eqref{thm:relative_perfectness_base_change2}.
	Moreover, \eqref{thm:relative_perfectness_base_change2} implies \eqref{thm:relative_perfectness_base_change3} by \cite[\href{https://stacks.math.columbia.edu/tag/01SI}{Tag 01SI}]{StacksProject}.
	Clearly, \eqref{thm:relative_perfectness_base_change3} implies \eqref{thm:relative_perfectness_base_change4}.
	By \Cref{lem:reflecting_relative_perfectness}\eqref{lem:reflecting_relative_perfectness3}, \eqref{thm:relative_perfectness_base_change4} implies \eqref{thm:relative_perfectness_base_change1}, and taking $t=1_S$ shows \eqref{thm:relative_perfectness_base_change1} implies \eqref{thm:relative_perfectness_base_change5}.

	We need to check \eqref{thm:relative_perfectness_base_change5} implies \eqref{thm:relative_perfectness_base_change4}.
	Let $t\colon T \to S$ be an affine surjective morphism such that $\mathbf{L}(t^\prime)^\ast K$ is relatively perfect over $Y_1\times_S T$.
	Choose a closed immersion $s\colon \operatorname{Spec}(k)\to S$ associated to a closed point of $S$.
	Consider the fibered square
	\begin{displaymath}
		\begin{tikzcd}
			{T\times_S \operatorname{Spec}(k)} & {\operatorname{Spec}(k)} \\
			T & {S.}
			\arrow["{u}", from=1-1, to=1-2]
			\arrow["{s^\prime}"', from=1-1, to=2-1]
			\arrow["s", from=1-2, to=2-2]
			\arrow["t"', from=2-1, to=2-2]
		\end{tikzcd}
	\end{displaymath}
	By base change, $s^\prime$ is a closed immersion and $u$ is affine and surjective.
	In particular, $T\times_S \operatorname{Spec}(k)$ is Noetherian, and $u$ is faithfully flat as its target is the spectrum of a field.
	Since $s^\prime$ is affine, \Cref{lem:relative_perfectness_affine_base_change} implies the derived pullback of $K$ to $Y_1\times_S Y_2\times_S T\times_S\operatorname{Spec}(k)$ is relatively perfect over $Y_1\times_S T\times_S \operatorname{Spec}(k)$.
	Since $u$ is affine and faithfully flat, \Cref{lem:reflecting_relative_perfectness}\eqref{lem:reflecting_relative_perfectness1}, applied to the families $Y_i\times_S \operatorname{Spec}(k)\to\operatorname{Spec}(k)$, ensures that the derived pullback of $K$ to $Y_1\times_S Y_2\times_S \operatorname{Spec}(k)$ is relatively perfect over $Y_1\times_S \operatorname{Spec}(k)$.
	As $s$ was an arbitrary closed immersion associated to a closed point, this proves \eqref{thm:relative_perfectness_base_change4}.
\end{proof}

\section{Mates}
\label{sec:mates}

\begin{remark}
	We recall the calculus of mates to the extent needed below. 
	See \cite[IV.7]{MacLane:1978}, \cite[\S 4]{Riehl:26}, and \cite[Sect.\ 1]{Cheng/Gurski/Riehl:2014} for details.
\end{remark}

\begin{setup}
	\label{setup:mates}
	Consider a diagram of functors and adjoint functors
	\begin{equation}
		\label{diag:mates}
		\begin{tikzcd}
			{\mathcal{C}} & {\mathcal{C}^\prime} \\
			{\mathcal{D}} & {\mathcal{D}^\prime,}
			\arrow["H", from=1-1, to=1-2]
			\arrow[""{name=0, anchor=center, inner sep=0}, "F"', shift right=2, from=1-1, to=2-1]
			\arrow[""{name=1, anchor=center, inner sep=0}, "{F^\prime}"', shift right=2, from=1-2, to=2-2]
			\arrow[""{name=2, anchor=center, inner sep=0}, "G"', shift right=2, from=2-1, to=1-1]
			\arrow["K"', from=2-1, to=2-2]
			\arrow[""{name=3, anchor=center, inner sep=0}, "{G^\prime}"', shift right=2, from=2-2, to=1-2]
			\arrow["\dashv"{anchor=center}, draw=none, from=0, to=2]
			\arrow["\dashv"{anchor=center}, draw=none, from=1, to=3]
		\end{tikzcd}
	\end{equation}
	where $F\dashv G$ has unit $\eta$ and counit $\epsilon$, and $F^\prime\dashv G^\prime$ has unit $\eta^\prime$ and counit $\epsilon^\prime$.
\end{setup}

\begin{definition}
	\label{def:mates}
	Assume \Cref{setup:mates}.
	Natural transformations $\alpha\colon F^\prime H\Rightarrow KF$ and $\beta\colon HG\Rightarrow G^\prime K$ are \textbf{mates} under \eqref{diag:mates} if $\beta$ equals the composite
	\begin{displaymath}
		HG\xrightarrow{\eta^\prime HG} G^\prime F^\prime HG \xrightarrow{G^\prime\alpha G} G^\prime KFG \xrightarrow{G^\prime K\epsilon} G^\prime K.
	\end{displaymath}
	Equivalently, $\alpha$ equals the composite
	\begin{displaymath}
		F^\prime H\xrightarrow{F^\prime H\eta} F^\prime HGF \xrightarrow{F^\prime\beta F} F^\prime G^\prime KF \xrightarrow{\epsilon^\prime KF} KF,
	\end{displaymath}
	see \cite[Sect.\ 1]{Cheng/Gurski/Riehl:2014}.
	Here we say $\alpha$ is the \textbf{left mate} of $\beta$ and $\beta$ is the \textbf{right mate} of $\alpha$.
\end{definition}

\begin{lemma}
	\label{lem:mates_equivalent_conditions}
	Assume \Cref{setup:mates}.
	Let $\alpha\colon F^\prime H\Rightarrow KF$ and $\beta\colon HG\Rightarrow G^\prime K$ be natural transformations.
	Then the following are equivalent:
	\begin{enumerate}
		\item \label{lem:mates_equivalent_conditions1} $\alpha$ and $\beta$ are mates
		\item \label{lem:mates_equivalent_conditions2} the following diagram commutes:
		\begin{displaymath}
			\begin{tikzcd}
				{F^\prime HG} & {KFG} \\
				{F^\prime G^\prime K} & {K}
				\arrow["{\alpha G}", from=1-1, to=1-2]
				\arrow["{F^\prime\beta}"', from=1-1, to=2-1]
				\arrow["{K\epsilon}", from=1-2, to=2-2]
				\arrow["{\epsilon^\prime K}"', from=2-1, to=2-2]
			\end{tikzcd}
		\end{displaymath}
		\item \label{lem:mates_equivalent_conditions3} the following diagram commutes:
		\begin{displaymath}
			\begin{tikzcd}
				H & {G^\prime F^\prime H} \\
				{HGF} & {G^\prime KF}
				\arrow["{\eta^\prime H}", from=1-1, to=1-2]
				\arrow["{H\eta}"', from=1-1, to=2-1]
				\arrow["{G^\prime\alpha}", from=1-2, to=2-2]
				\arrow["{\beta F}"', from=2-1, to=2-2]
			\end{tikzcd}
		\end{displaymath}
		\item \label{lem:mates_equivalent_conditions4} for all objects $c\in\mathcal{C}$ and $d\in\mathcal{D}$, the following diagram commutes:
		\begin{displaymath}
			\begin{tikzcd}
				{\operatorname{Hom}_{\mathcal{D}}(Fc,d)} & {\operatorname{Hom}_{\mathcal{C}}(c,Gd)} \\
				{\operatorname{Hom}_{\mathcal{D}^\prime}(KFc,Kd)} & {\operatorname{Hom}_{\mathcal{C}^\prime}(Hc,HGd)} \\
				{\operatorname{Hom}_{\mathcal{D}^\prime}(F^\prime Hc,Kd)} & {\operatorname{Hom}_{\mathcal{C}^\prime}(Hc,G^\prime Kd)}
				\arrow["\cong", from=1-1, to=1-2]
				\arrow["K"', from=1-1, to=2-1]
				\arrow["H", from=1-2, to=2-2]
				\arrow["{-\circ\alpha_c}"', from=2-1, to=3-1]
				\arrow["{\beta_d\circ -}", from=2-2, to=3-2]
				\arrow["\cong", from=3-1, to=3-2]
			\end{tikzcd}
		\end{displaymath}
		where the top and bottom bijections are induced by the adjunctions $F\dashv G$ and $F^\prime\dashv G^\prime$, respectively.
	\end{enumerate}
\end{lemma}

\begin{proof}
	The equivalence of \eqref{lem:mates_equivalent_conditions3} and \eqref{lem:mates_equivalent_conditions4} is \cite[IV.7, Exercise 4]{MacLane:1978}; we recall the argument.
	Let $c\in\mathcal{C}$, $d\in\mathcal{D}$, and $\phi\colon Fc\to d$.
	Going down the left column and then along the bottom row of the diagram in \eqref{lem:mates_equivalent_conditions4}, $\phi$ is sent to $G^\prime K(\phi)\circ G^\prime(\alpha_c)\circ\eta^\prime_{Hc}$.
	Going along the top row and then down the right column, $\phi$ is sent to
	\begin{displaymath}
		\beta_d\circ HG(\phi)\circ H(\eta_c) = G^\prime K(\phi)\circ\beta_{Fc}\circ H(\eta_c),
	\end{displaymath}
	where the equality is naturality of $\beta$.
	Hence, \eqref{lem:mates_equivalent_conditions3} implies \eqref{lem:mates_equivalent_conditions4}.
	Conversely, taking $d=Fc$ and $\phi=1_{Fc}$ shows that \eqref{lem:mates_equivalent_conditions4} implies \eqref{lem:mates_equivalent_conditions3}.

	Next, the right mate of $\alpha$ equals $\beta$ if, and only if, the corresponding adjunct morphisms $F^\prime HG\Rightarrow K$ with respect to $F^\prime\dashv G^\prime$ are equal, which is exactly \eqref{lem:mates_equivalent_conditions2}.
	Dually, the left mate of $\beta$ equals $\alpha$ if, and only if, the corresponding adjunct morphisms $H\Rightarrow G^\prime KF$ with respect to $F^\prime\dashv G^\prime$ are equal, which is exactly \eqref{lem:mates_equivalent_conditions3}.
	This shows that \eqref{lem:mates_equivalent_conditions1}, \eqref{lem:mates_equivalent_conditions2}, and \eqref{lem:mates_equivalent_conditions3} are equivalent.
\end{proof}

\begin{remark}
	\label{rem:mates_of_isomorphisms}
	The mate of a natural isomorphism need not be an isomorphism.
	For instance, the unit of an adjunction is the right mate of the identity transformation of the left adjoint (take $H$ to be an identity functor and the adjunction on the left of \eqref{diag:mates} to be the identity adjunction).
	Hence, any adjunction whose unit is not an isomorphism (e.g.\ the free-forgetful adjunction between sets and groups) gives an example.
\end{remark}

\section{Base change for adjoints}
\label{sec:base_change_adjoints}

\begin{lemma}
	\label{lem:internal_hom_pullback}
	Let $f\colon W \to Y$ be a morphism of schemes.
	Let $E,G\in D_{\operatorname{qc}}(Y)$ where $E$ is perfect.
	Then there exists an isomorphism
	\begin{displaymath}
		\mathbf{L}f^\ast \operatorname{\mathbf{R}\mathcal{H}\! \mathit{om}} (E,G) \to  \operatorname{\mathbf{R}\mathcal{H}\! \mathit{om}} (\mathbf{L}f^\ast E , \mathbf{L}f^\ast G).
	\end{displaymath}
\end{lemma}

\begin{proof}
	By \cite[\href{https://stacks.math.columbia.edu/tag/0A6H}{Tag 0A6H}]{StacksProject}, there are isomorphisms
	\begin{displaymath}
		\operatorname{\mathbf{R}\mathcal{H}\! \mathit{om}}(E,G)
		\cong \operatorname{\mathbb{R}\mathcal{H}\! \mathit{om}}(E,G)
		\quad\textrm{and}\quad
		\operatorname{\mathbf{R}\mathcal{H}\! \mathit{om}}(\mathbf{L}f^\ast E, \mathbf{L}f^\ast G)
		\cong \operatorname{\mathbb{R}\mathcal{H}\! \mathit{om}}(\mathbf{L}f^\ast E,\mathbf{L}f^\ast G).
	\end{displaymath}
	By \cite[\href{https://stacks.math.columbia.edu/tag/08I3}{Tag 08I3}]{StacksProject}, there exists a natural morphism
	\begin{displaymath}
		\mathbf{L}f^\ast \operatorname{\mathbb{R}\mathcal{H}\! \mathit{om}} (E,G) \to  \operatorname{\mathbb{R}\mathcal{H}\! \mathit{om}} (\mathbf{L}f^\ast E , \mathbf{L}f^\ast G).
	\end{displaymath}
	The claim follows if we can show this is an isomorphism, which is a special case of \cite[Proposition 22.70]{Gortz/Wedhorn:2023}.
\end{proof}

\begin{lemma}
	[Base change for relative duals]
	\label{lem:base_change_relative_dual}
	Consider a fibered square of Noetherian schemes
	\begin{displaymath}
		\begin{tikzcd}
			{X^\prime} & {S^\prime} \\
			X & S
			\arrow["{f^\prime}", from=1-1, to=1-2]
			\arrow["{g^\prime}"', from=1-1, to=2-1]
			\arrow["g", from=1-2, to=2-2]
			\arrow["f"', from=2-1, to=2-2]
		\end{tikzcd}
	\end{displaymath}
	where $f$ is proper and flat.
	Let
	\begin{displaymath}
		\rho_{f,g}\colon \mathbf{L}(g^\prime)^\ast f^!\mathcal{O}_S
		\to (f^\prime)^!\mathcal{O}_{S^\prime}
	\end{displaymath}
	be the canonical base change morphism.
	Then $\rho_{f,g}$ is an isomorphism and $\mathbf{L}(g^\prime)^\ast f^!\mathcal{O}_S$ is a relative dualizing complex for $f^\prime$.
	Furthermore, if $g$ is affine and $E\in D^b_{\operatorname{coh}}(X)$ is $f$-perfect, then the canonical morphism
	\begin{equation}
		\label{eq:base_change_relative_dual}
		\begin{aligned}
			\theta_E
			& \colon \mathbf{L}(g^\prime)^\ast
			\operatorname{\mathbf{R}\mathcal{H}\! \mathit{om}}
			(E,f^!\mathcal{O}_S)
			\to \operatorname{\mathbf{R}\mathcal{H}\! \mathit{om}}
			(\mathbf{L}(g^\prime)^\ast E,\mathbf{L}(g^\prime)^\ast f^!\mathcal{O}_S)
			\\& \xrightarrow{ \operatorname{\mathbf{R}\mathcal{H}\! \mathit{om}}(1,\rho_{f,g})} \operatorname{\mathbf{R}\mathcal{H}\! \mathit{om}}
			(\mathbf{L}(g^\prime)^\ast E, (f^\prime)^!\mathcal{O}_{S^\prime})
		\end{aligned}
	\end{equation}
	is an isomorphism, where the first morphism in \eqref{eq:base_change_relative_dual} is the canonical morphism of \cite[\href{https://stacks.math.columbia.edu/tag/08I3}{Tag 08I3}]{StacksProject}.
\end{lemma}

\begin{proof}
	Since $f$ and $f^\prime$ are proper, we may identify $f^\times\mathcal{O}_S$ with $f^!\mathcal{O}_S$ and $(f^\prime)^\times\mathcal{O}_{S^\prime}$ with $(f^\prime)^!\mathcal{O}_{S^\prime}$ \cite[\href{https://stacks.math.columbia.edu/tag/0AU3}{Tag 0AU3}]{StacksProject}.
	Then the first claim is \cite[\href{https://stacks.math.columbia.edu/tag/0AAB}{Tag 0AAB} and \href{https://stacks.math.columbia.edu/tag/0B6S}{0B6S}]{StacksProject}.
	See also \cite[\href{https://stacks.math.columbia.edu/tag/0E2Y}{Tag 0E2Y} and \href{https://stacks.math.columbia.edu/tag/0E2Z}{0E2Z}]{StacksProject} for relative dualizing complexes.
	Moreover, $\rho_{f,g}$ is compatible with the trace morphisms for $f$ and $f^\prime$ \cite[\href{https://stacks.math.columbia.edu/tag/0B6J}{Tag 0B6J}]{StacksProject}.

	We prove the second claim.
	The complexes occurring in \eqref{eq:base_change_relative_dual} belong to $D_{\operatorname{qc}}(X^\prime)$.
	Indeed, derived pullback preserves pseudocoherence \cite[\href{https://stacks.math.columbia.edu/tag/08H4}{Tag 08H4}]{StacksProject}, and upper shriek for a proper morphism preserves complexes with bounded below and coherent cohomology \cite[\href{https://stacks.math.columbia.edu/tag/0AU1}{Tag 0AU1}]{StacksProject}.
	Then apply \cite[\href{https://stacks.math.columbia.edu/tag/0A6H}{Tag 0A6H}]{StacksProject}.

	Choose a compact generator $P$ of $D_{\operatorname{qc}}(X)$ \cite[\href{https://stacks.math.columbia.edu/tag/09IS}{Tag 09IS}]{StacksProject}.
	Since $g^\prime$ is affine, $\mathbf{L}(g^\prime)^\ast P$ is a compact generator of $D_{\operatorname{qc}}(X^\prime)$ \cite[\href{https://stacks.math.columbia.edu/tag/0BQT}{Tag 0BQT}]{StacksProject}.
	Thus, it is enough to prove
	\begin{equation}
		\label{eq:pushforward_generator_test}
		\mathbf{R}(f^\prime)_\ast ( \operatorname{\mathbf{R}\mathcal{H}\! \mathit{om}} (\mathbf{L}(g^\prime)^\ast P,\mathcal{O}_{X^\prime}) \otimes^{\mathbf{L}}\operatorname{cone}(\theta_E))
		\cong 0.
	\end{equation}
	Indeed, if \eqref{eq:pushforward_generator_test} holds, then for every $n\in\mathbb{Z}$,
	\begin{displaymath}
		\begin{aligned}
			\operatorname{Hom}
			&(\mathbf{L}(g^\prime)^\ast P,\operatorname{cone}(\theta_E)[n])
			\\&\cong \operatorname{Hom}( \mathcal{O}_{X^\prime},\operatorname{\mathbf{R}\mathcal{H}\! \mathit{om}} (\mathbf{L}(g^\prime)^\ast P,\mathcal{O}_{X^\prime}) \otimes^{\mathbf{L}}\operatorname{cone}(\theta_E)[n])
			\\&\cong\operatorname{Hom}( \mathcal{O}_{S^\prime}, \mathbf{R}(f^\prime)_\ast( \operatorname{\mathbf{R}\mathcal{H}\! \mathit{om}} (\mathbf{L}(g^\prime)^\ast P,\mathcal{O}_{X^\prime})\otimes^{\mathbf{L}}\operatorname{cone}(\theta_E))[n])
			\\&\cong 0.
		\end{aligned}
	\end{displaymath}
	Since $\mathbf{L}(g^\prime)^\ast P$ generates $D_{\operatorname{qc}}(X^\prime)$, it follows that $\operatorname{cone}(\theta_E)\cong 0$.
	Moreover, by \cite[\href{https://stacks.math.columbia.edu/tag/08DJ}{Tag 08DJ} and \href{https://stacks.math.columbia.edu/tag/08DQ}{08DQ}]{StacksProject} and \Cref{lem:internal_hom_pullback}, $\theta_{P\otimes^{\mathbf{L}}E}$ identifies with the morphism
	\begin{displaymath}
		\operatorname{\mathbf{R}\mathcal{H}\! \mathit{om}}
		(\mathbf{L}(g^\prime)^\ast P,\mathcal{O}_{X^\prime})
		\otimes^{\mathbf{L}}\theta_E.
	\end{displaymath}
	Hence, it suffices to show that $\mathbf{R}(f^\prime)_\ast(\theta_{P\otimes^{\mathbf{L}}E})$ is an isomorphism.

	Since $E$ is $f$-perfect and $P$ is perfect, $\mathbf{R}f_\ast(P\otimes^{\mathbf{L}}E)\in\operatorname{Perf}(S)$.
	For any $B\in D^-_{\operatorname{coh}}(X)$, denote by
	\begin{displaymath}
		\xi_f(B)\colon
		\mathbf{R}f_\ast \operatorname{\mathbf{R}\mathcal{H}\! \mathit{om}}  (B,f^!\mathcal{O}_S)
		\to \operatorname{\mathbf{R}\mathcal{H}\! \mathit{om}}  (\mathbf{R}f_\ast B,\mathcal{O}_S)
	\end{displaymath}
	the duality morphism from \cite[\href{https://stacks.math.columbia.edu/tag/0A9Q}{Tag 0A9Q} and \href{https://stacks.math.columbia.edu/tag/0GEV}{0GEV}]{StacksProject}.
	Define $\xi_{f^\prime}$ similarly.
	By \cite[\href{https://stacks.math.columbia.edu/tag/0GEW}{Tag 0GEW}]{StacksProject}, these are isomorphisms.
	By the construction of $\xi_f(P\otimes^{\mathbf{L}}E)$ \cite[\href{https://stacks.math.columbia.edu/tag/0B69}{Tag 0B69}]{StacksProject}, its adjunct is
	\begin{equation}
		\label{eq:duality_via_trace}
		\begin{aligned}
			&\mathbf{R}f_\ast \operatorname{\mathbf{R}\mathcal{H}\! \mathit{om}} (P\otimes^{\mathbf{L}}E,f^!\mathcal{O}_S) \otimes^{\mathbf{L}} \mathbf{R}f_\ast(P\otimes^{\mathbf{L}}E)
			\\&\to \mathbf{R}f_\ast(\operatorname{\mathbf{R}\mathcal{H}\! \mathit{om}}(P\otimes^{\mathbf{L}}E,f^!\mathcal{O}_S) \otimes^{\mathbf{L}}(P\otimes^{\mathbf{L}}E))
			\\&\to\mathbf{R}f_\ast f^!\mathcal{O}_S \xrightarrow{\operatorname{Tr}_f} \mathcal{O}_S.
		\end{aligned}
	\end{equation}
	The analogous description holds for $\xi_{f^\prime}(\mathbf{L}(g^\prime)^\ast(P\otimes^{\mathbf{L}}E))$.

	For $B\in D_{\operatorname{qc}}(X)$, denote by
	\begin{displaymath}
		\beta_B\colon
		\mathbf{L}g^\ast\mathbf{R}f_\ast B
		\to \mathbf{R}(f^\prime)_\ast\mathbf{L}(g^\prime)^\ast B
	\end{displaymath}
	the canonical base change morphism \cite[\href{https://stacks.math.columbia.edu/tag/08HY}{Tag 08HY}]{StacksProject}.
	As $f$ is flat, $\beta_B$ is an isomorphism \cite[\href{https://stacks.math.columbia.edu/tag/08IB}{Tag 08IB}]{StacksProject}.
	We claim that the following diagram commutes:
	\begin{equation}
		\label{diag:relative_dual_base_change_compatibility}
		\begin{tikzcd}[column sep=small]
			\mathbf{L}g^\ast\mathbf{R}f_\ast
			\operatorname{\mathbf{R}\mathcal{H}\! \mathit{om}}
			(P\otimes^{\mathbf{L}}E,f^!\mathcal{O}_S)
			\arrow[r, "{\mathbf{L}g^\ast \xi_f(P\otimes^{\mathbf{L}}E)}"]
			\arrow[d, "{\beta_{\operatorname{\mathbf{R}\mathcal{H}\! \mathit{om}}(P\otimes^{\mathbf{L}}E,f^!\mathcal{O}_S)}}"']
			&\mathbf{L}g^\ast \operatorname{\mathbf{R}\mathcal{H}\! \mathit{om}}
			( \mathbf{R}f_\ast(P\otimes^{\mathbf{L}}E), \mathcal{O}_S)
			\arrow[d]
			\\\mathbf{R}(f^\prime)_\ast\mathbf{L}(g^\prime)^\ast \operatorname{\mathbf{R}\mathcal{H}\! \mathit{om}}(P\otimes^{\mathbf{L}}E,f^!\mathcal{O}_S)
			\arrow[d, "{\mathbf{R}(f^\prime)_\ast(\theta_{P\otimes^{\mathbf{L}}E})}"']
			&\operatorname{\mathbf{R}\mathcal{H}\! \mathit{om}}( \mathbf{L}g^\ast \mathbf{R}f_\ast(P\otimes^{\mathbf{L}}E), \mathcal{O}_{S^\prime})
			\arrow[d, "{\operatorname{\mathbf{R}\mathcal{H}\! \mathit{om}}(\beta_{P\otimes^{\mathbf{L}}E}^{-1},1)}"]
			\\\mathbf{R}(f^\prime)_\ast \operatorname{\mathbf{R}\mathcal{H}\! \mathit{om}}
			( \mathbf{L}(g^\prime)^\ast(P\otimes^{\mathbf{L}}E), (f^\prime)^!\mathcal{O}_{S^\prime})
			\arrow[r, "{\xi_{f^\prime} (\mathbf{L}(g^\prime)^\ast (P\otimes^{\mathbf{L}}E))}"' {yshift=-4pt}]
			&\operatorname{\mathbf{R}\mathcal{H}\! \mathit{om}}(  \mathbf{R}(f^\prime)_\ast \mathbf{L}(g^\prime)^\ast (P\otimes^{\mathbf{L}}E), \mathcal{O}_{S^\prime}).
		\end{tikzcd}
	\end{equation}
	The unlabeled right vertical morphism is the canonical morphism of \cite[\href{https://stacks.math.columbia.edu/tag/08I3}{Tag 08I3}]{StacksProject}.
	It is an isomorphism because $\mathbf{R}f_\ast(P\otimes^{\mathbf{L}}E)$ is perfect \cite[\href{https://stacks.math.columbia.edu/tag/08DQ}{Tag 08DQ}]{StacksProject}.

	We check the claimed compatibility.
	By \cite[\href{https://stacks.math.columbia.edu/tag/08DJ}{Tag 08DJ}]{StacksProject}, commutativity of \eqref{diag:relative_dual_base_change_compatibility} is equivalent to equality of the resulting two adjunct morphisms
	\begin{displaymath}
		\mathbf{L}g^\ast\mathbf{R}f_\ast \operatorname{\mathbf{R}\mathcal{H}\! \mathit{om}} (P\otimes^{\mathbf{L}}E,f^!\mathcal{O}_S) \otimes^{\mathbf{L}}  \mathbf{R}(f^\prime)_\ast \mathbf{L}(g^\prime)^\ast(P\otimes^{\mathbf{L}}E)
		\to \mathcal{O}_{S^\prime}.
	\end{displaymath}
	The morphism obtained from the upper right boundary is the composite
	\begin{displaymath}
		\begin{aligned}
			& \mathbf{L}g^\ast\mathbf{R}f_\ast \operatorname{\mathbf{R}\mathcal{H}\! \mathit{om}} (P\otimes^{\mathbf{L}}E,f^!\mathcal{O}_S) \otimes^{\mathbf{L}} \mathbf{R}(f^\prime)_\ast \mathbf{L}(g^\prime)^\ast(P\otimes^{\mathbf{L}}E)
			\\&\xrightarrow{1\otimes\beta_{P\otimes^{\mathbf{L}}E}^{-1}} \mathbf{L}g^\ast\mathbf{R}f_\ast \operatorname{\mathbf{R}\mathcal{H}\! \mathit{om}} (P\otimes^{\mathbf{L}}E,f^!\mathcal{O}_S) \otimes^{\mathbf{L}} \mathbf{L}g^\ast\mathbf{R}f_\ast(P\otimes^{\mathbf{L}}E)
			\\&\to \mathbf{L}g^\ast\mathbf{R}f_\ast f^!\mathcal{O}_S \xrightarrow{\mathbf{L}g^\ast\operatorname{Tr}_f} \mathcal{O}_{S^\prime},
		\end{aligned}
	\end{displaymath}
	where the unlabeled morphism is the pullback of the relative cup product followed by evaluation.
	The relative cup product is compatible with base change by \cite[\href{https://stacks.math.columbia.edu/tag/0H99}{Tag 0H99}]{StacksProject}.
	Moreover, compatibility of $\rho_{f,g}$ with the trace morphisms is expressed by the commutative diagram
	\begin{displaymath}
		\begin{tikzcd}[column sep=huge]
			\mathbf{L}g^\ast\mathbf{R}f_\ast f^!\mathcal{O}_S
			\arrow[r, "{\mathbf{L}g^\ast\operatorname{Tr}_f}"]
			\arrow[d,"{\beta_{f^!\mathcal{O}_S}}"']
			&\mathcal{O}_{S^\prime}
			\\{\mathbf{R}(f^\prime)_\ast \mathbf{L}(g^\prime)^\ast f^!\mathcal{O}_S}
			\arrow[r,"{\mathbf{R}(f^\prime)_\ast(\rho_{f,g})}"']
			& \mathbf{R}(f^\prime)_\ast (f^\prime)^!\mathcal{O}_{S^\prime}
			\arrow[u,"{\operatorname{Tr}_{f^\prime}}"']
		\end{tikzcd}
	\end{displaymath}
	by \cite[\href{https://stacks.math.columbia.edu/tag/0B6J}{Tag 0B6J}]{StacksProject}.
	On the other hand, the first arrow in $\theta_{P\otimes^{\mathbf{L}}E}$ arises from \cite[\href{https://stacks.math.columbia.edu/tag/08I3}{Tag 08I3}]{StacksProject}, which is constructed by pulling back the evaluation morphism.
	Its compatibility with the base change morphisms above is precisely the construction in \cite[\href{https://stacks.math.columbia.edu/tag/08I4}{Tag 08I4}]{StacksProject}.
	The second arrow in $\theta_{P\otimes^{\mathbf{L}}E}$ is induced by $\rho_{f,g}$.
	Together with the compatibility of relative cup products and traces above, the two adjunct morphisms coincide.
	Hence, \eqref{diag:relative_dual_base_change_compatibility} commutes.

	In \eqref{diag:relative_dual_base_change_compatibility}, every morphism, except possibly the lower left vertical one, is an isomorphism.
	Thus, $\mathbf{R}(f^\prime)_\ast(\theta_{P\otimes^{\mathbf{L}}E})$ is an isomorphism.
	Consequently,
	\begin{displaymath}
		\mathbf{R}(f^\prime)_\ast ( \operatorname{\mathbf{R}\mathcal{H}\! \mathit{om}} (\mathbf{L}(g^\prime)^\ast P,\mathcal{O}_{X^\prime}) \otimes^{\mathbf{L}}\theta_E )
	\end{displaymath}
	is an isomorphism, which completes the proof.
\end{proof}

\begin{lemma}
	\label{lem:relative_dual_relatively_perfect}
	Let $f\colon Y\to X$ be a proper morphism of Noetherian schemes where $X$ is separated.
	If $E\in D^-_{\operatorname{coh}}(Y)$ is $f$-perfect, then $\operatorname{\mathbf{R}\mathcal{H}\! \mathit{om}}(E, f^! P)$ is $f$-perfect for all $P\in \operatorname{Perf}(X)$.
\end{lemma}

\begin{proof}
	As $f$ is proper and $X$ is separated, $Y$ is separated.
	Hence, the claim is \cite[Lemma 3.7]{Ballard:2009}.
\end{proof}

\begin{proposition}
	\label{prop:right_adjoint_base_change}
	Consider \Cref{setup:fm_cube} with $t$ affine.
	Let $K\in D^b_{\operatorname{coh}}(Y_1\times_S Y_2)$ be relatively perfect over $Y_2$.
	Then $\Phi_{\mathbf{L}(t^\prime)^\ast K^\prime}$ is right adjoint to $\Phi_{\mathbf{L} (t^\prime)^\ast K}$ on $D_{\operatorname{qc}}$ where
	\begin{displaymath}
		K^\prime:= \operatorname{\mathbf{R}\mathcal{H}\! \mathit{om}} (  K,  (f^\prime_1)^! \mathcal{O}_{Y_2}).
	\end{displaymath}
	Moreover, $K^\prime\in D^b_{\operatorname{coh}}(Y_1\times_S Y_2)$ and $\mathbf{L}(t^\prime)^\ast K^\prime\in D^b_{\operatorname{coh}}(Y_1\times_S Y_2\times_S T)$.
\end{proposition}

\begin{proof}
	We check the first claim.
	For any $E\in D_{\operatorname{qc}}(Y_1\times_S T)$ and $A\in D_{\operatorname{qc}}(Y_2\times_S T)$, there is a string of natural isomorphisms obtained by adjunctions:
	\begin{displaymath}
		\begin{aligned}
			\operatorname{Hom} & ( \Phi_{\mathbf{L} (t^\prime)^\ast K} (E) , A )
			\\&= \operatorname{Hom} ( \mathbf{R} (g^\prime_1)_\ast (\mathbf{L} (g^\prime_2)^\ast E \otimes^{\mathbf{L}} \mathbf{L} (t^\prime)^\ast K), A ) && \textrm{(definition)}
			\\&\cong \operatorname{Hom} ( \mathbf{L} (g^\prime_2)^\ast E \otimes^{\mathbf{L}} \mathbf{L} (t^\prime)^\ast K, (g^\prime_1)^\times A ) && \textrm{(\cite[\href{https://stacks.math.columbia.edu/tag/0A9E}{Tag 0A9E}]{StacksProject})}
			\\&\cong \operatorname{Hom} ( \mathbf{L} (g^\prime_2)^\ast  E , \operatorname{\mathbf{R}\mathcal{H}\! \mathit{om}} ( \mathbf{L} (t^\prime)^\ast K, (g^\prime_1)^\times A)) && \textrm{(\cite[\href{https://stacks.math.columbia.edu/tag/08DH}{Tag 08DH}]{StacksProject})}
			\\&\cong \operatorname{Hom} ( E , \mathbf{R} (g^\prime_2)_\ast \operatorname{\mathbf{R}\mathcal{H}\! \mathit{om}} ( \mathbf{L} (t^\prime)^\ast K, (g^\prime_1)^\times A)) && \textrm{(adjunction)}.
		\end{aligned}
	\end{displaymath}
	Here $(g_1^\prime)^\times$ denotes the right adjoint of $\mathbf{R}(g_1^\prime)_\ast$ on $D_{\operatorname{qc}}$.
	The desired claim follows if we can find an isomorphism
	\begin{displaymath}
		\begin{aligned}
			\mathbf{R} & (g^\prime_2)_\ast \operatorname{\mathbf{R}\mathcal{H}\! \mathit{om}} ( \mathbf{L} (t^\prime)^\ast K, (g^\prime_1)^\times A) \\&\cong \mathbf{R} (g^\prime_2)_\ast ( \mathbf{L} (t^\prime)^\ast \operatorname{\mathbf{R}\mathcal{H}\! \mathit{om}} (  K,  (f^\prime_1)^! \mathcal{O}_{Y_2}) \otimes^\mathbf{L} \mathbf{L} (g^\prime_1)^\ast A ) = \Phi_{\mathbf{L}(t^\prime)^\ast K^\prime} (A).
		\end{aligned}
	\end{displaymath}
	Moreover, $(f^\prime_1)^! \mathcal{O}_{Y_2} \in D^+_{\operatorname{coh}}(Y_1\times_S Y_2)$ as upper shriek functors for proper morphisms preserve complexes with bounded below and coherent cohomology \cite[\href{https://stacks.math.columbia.edu/tag/0AU1}{Tag 0AU1}]{StacksProject}.
	By \Cref{lem:base_change_relative_dual},
	\begin{equation}
		\label{eq:right_adjoint_kernel_base_change}
		\mathbf{L}(t^\prime)^\ast
		\operatorname{\mathbf{R}\mathcal{H}\! \mathit{om}}( K, (f^\prime_1)^! \mathcal{O}_{Y_2})
		\cong \operatorname{\mathbf{R}\mathcal{H}\! \mathit{om}}(\mathbf{L}(t^\prime)^\ast K, (g^\prime_1)^! \mathcal{O}_{Y_2\times_S T}).
	\end{equation}
	By \Cref{lem:relative_perfectness_affine_base_change}, $\mathbf{L}(t^\prime)^\ast K$ is relatively perfect over $Y_2\times_S T$, i.e.\
	\begin{displaymath}
		\mathbf{R}(g_1^\prime)_\ast (\mathbf{L}(t^\prime)^\ast K \otimes^{\mathbf{L}} \operatorname{Perf}(Y_1\times_S Y_2 \times_S T))\subseteq \operatorname{Perf}(Y_2\times_S T).
	\end{displaymath}
	There is another string of isomorphisms for all $A \in D_{\operatorname{qc}}(Y_2\times_S T)$:
	\begin{displaymath}
		\begin{aligned}
			\mathbf{R} & (g^\prime_2)_\ast \operatorname{\mathbf{R}\mathcal{H}\! \mathit{om}} ( \mathbf{L} (t^\prime)^\ast K, (g^\prime_1)^\times A)
			\\&\cong \mathbf{R} (g^\prime_2)_\ast \operatorname{\mathbf{R}\mathcal{H}\! \mathit{om}} ( \mathbf{L} (t^\prime)^\ast K, \mathbf{L} (g^\prime_1)^\ast A \otimes^\mathbf{L} (g^\prime_1)^! \mathcal{O}_{Y_2 \times_S T} ) && \textrm{(\cite[\href{https://stacks.math.columbia.edu/tag/0B6S}{Tag 0B6S}]{StacksProject})}
			\\&\cong \mathbf{R} (g^\prime_2)_\ast ( \operatorname{\mathbf{R}\mathcal{H}\! \mathit{om}} ( \mathbf{L} (t^\prime)^\ast K, (g^\prime_1)^! \mathcal{O}_{Y_2 \times_S T} ) \otimes^\mathbf{L} \mathbf{L} (g^\prime_1)^\ast A ) && \textrm{(\cite[Lemma 2.13]{Rizzardo:2017})}
			\\&\cong \mathbf{R} (g^\prime_2)_\ast ( \mathbf{L} (t^\prime)^\ast \operatorname{\mathbf{R}\mathcal{H}\! \mathit{om}} (  K,  (f^\prime_1)^! \mathcal{O}_{Y_2}) \otimes^\mathbf{L} \mathbf{L} (g^\prime_1)^\ast A ) && \textrm{(use \eqref{eq:right_adjoint_kernel_base_change})}.
		\end{aligned}
	\end{displaymath}
	The first isomorphism also uses that $(g_1^\prime)^\times\mathcal{O}_{Y_2\times_S T}\cong (g_1^\prime)^!\mathcal{O}_{Y_2\times_S T}$ \cite[\href{https://stacks.math.columbia.edu/tag/0AU3}{Tag 0AU3}]{StacksProject}, and the second uses the fact that $\mathbf{L}(t^\prime)^\ast K$ is $g_1^\prime$-perfect.
	This proves the first claim.

	Lastly, \cite[Definition 2.11 and Lemma 2.12]{Rizzardo:2017} imply $K^\prime\in D^b_{\operatorname{coh}}(Y_1\times_S Y_2)$.
	Applying loc.\ cit.\ to $\mathbf{L}(t^\prime)^\ast K$, which is relatively perfect over $Y_2\times_S T$, the right-hand side of \eqref{eq:right_adjoint_kernel_base_change} belongs to $D^b_{\operatorname{coh}}(Y_1\times_S Y_2\times_S T)$, and hence so does $\mathbf{L}(t^\prime)^\ast K^\prime$.
\end{proof}

\begin{corollary}
	\label{cor:right_adjoint_restricts}
	Consider \Cref{setup:fm_cube} with $S$ separated and $t=1_S$.
	Let $K\in D^b_{\operatorname{coh}}(Y_1\times_S Y_2)$ be relatively perfect over $Y_2$.
	Denote by $K^\prime$ the kernel of the integral transform obtained in \Cref{prop:right_adjoint_base_change} which is right adjoint to $\Phi_K$ on $D_{\operatorname{qc}}$.
	Then $K^\prime$ is relatively perfect over $Y_2$, and so $\Phi_{K^\prime}$ restricts to $D^b_{\operatorname{coh}}$.
	In particular, if $K$ is relatively perfect over each $Y_i$, then $\Phi_K$ and $\Phi_{K^\prime}$ restrict to an adjoint pair on $D^b_{\operatorname{coh}}$.
\end{corollary}

\begin{proof}
	As $Y_2$ is separated (\Cref{rem:fm_cube_separated}), \Cref{lem:relative_dual_relatively_perfect} applies to $f_1^\prime$ and $\mathcal{O}_{Y_2}$, and so the kernel $K^\prime$ is $f_1^\prime$-perfect, i.e.\ relatively perfect over $Y_2$.
	Since $Y_2$ is the source of the integral transform $\Phi_{K^\prime}\colon D_{\operatorname{qc}}(Y_2)\to D_{\operatorname{qc}}(Y_1)$, it follows from \cite[Propositions 3.2 and 3.6]{Dutta/Lank/ManaliRahul:2025} that $\Phi_{K^\prime}(D^b_{\operatorname{coh}}(Y_2))\subseteq D^b_{\operatorname{coh}}(Y_1)$.
	Moreover, if $K$ is relatively perfect over $Y_1$, then $\Phi_K(D^b_{\operatorname{coh}}(Y_1))\subseteq D^b_{\operatorname{coh}}(Y_2)$, and the adjunction on $D_{\operatorname{qc}}$ restricts to the full subcategories $D^b_{\operatorname{coh}}$.
\end{proof}

\begin{lemma}
	\label{lem:left_adjoint}
	Consider \Cref{setup:fm_cube} with $S$ separated.
	Suppose $K\in D^b_{\operatorname{coh}}(Y_1\times_S Y_2)$ is relatively perfect over $Y_1$.
	Then $\Phi_K$ admits a left adjoint on $D_{\operatorname{qc}}$.
	In particular, the left adjoint is the integral transform with kernel $\operatorname{\mathbf{R}\mathcal{H}\! \mathit{om}} (K, (f^\prime_2)^! \mathcal{O}_{Y_1})$.
\end{lemma}

\begin{proof}
	This is in the proof of \cite[Proposition 4.3]{Dutta/Lank/ManaliRahul:2025}.
	We spell it out for convenience.
	By \Cref{rem:fm_cube_separated}, $Y_1$ and $Y_1\times_S Y_2$ are separated, so \cite[Lemma 3.10]{Ballard:2009} applies to $f_2^\prime$.
	There exists a string of natural isomorphisms for all $E\in D_{\operatorname{qc}}(Y_2)$ and $G\in D_{\operatorname{qc}}(Y_1)$:
	\begin{displaymath}
		\begin{aligned}
			\operatorname{Hom}& (E,\Phi_K (G))
			\cong \operatorname{Hom}( E ,\mathbf{R}(f_1^\prime)_\ast (K \otimes^{\mathbf{L}} \mathbf{L}(f^\prime_2)^\ast G)) && \textrm{(definition)}
			\\&\cong \operatorname{Hom}( \mathbf{L}(f_1^\prime)^\ast E , K \otimes^{\mathbf{L}} \mathbf{L}(f_2^\prime)^\ast G) && \textrm{(adjunction)}
			\\&\cong \operatorname{Hom}( \mathbf{R}(f_2^\prime)_\ast ( \operatorname{\mathbf{R}\mathcal{H}\! \mathit{om}} (K, (f^\prime_2)^! \mathcal{O}_{Y_1}) \otimes^{\mathbf{L}} \mathbf{L}(f_1^\prime)^\ast E ) , G) && \textrm{(\cite[Lemma 3.10]{Ballard:2009})}.
		\end{aligned}
	\end{displaymath}
	Thus, the claim follows.
\end{proof}

\begin{proposition}
	\label{prop:left_adjoint_base_change}
	Consider \Cref{setup:fm_cube} where $S$ is separated and $t$ is affine.
	Let $K\in D^b_{\operatorname{coh}}(Y_1\times_S Y_2)$ be relatively perfect over $Y_1$.
	Then $\Phi_{\mathbf{L}(t^\prime)^\ast K^{\prime\prime}}$ is left adjoint to $\Phi_{\mathbf{L} (t^\prime)^\ast K}$ on $D_{\operatorname{qc}}$ where
	\begin{displaymath}
		K^{\prime\prime}:= \operatorname{\mathbf{R}\mathcal{H}\! \mathit{om}} (  K,  (f^\prime_2)^! \mathcal{O}_{Y_1}).
	\end{displaymath}
	In such a case, $K^{\prime\prime}\in D^b_{\operatorname{coh}}(Y_1\times_S Y_2)$.
\end{proposition}

\begin{proof}
	The last claim follows from \cite[Definition 2.11 and Lemma 2.12]{Rizzardo:2017}.
	We check the first claim.
	By \Cref{lem:relative_perfectness_affine_base_change}, $\mathbf{L} (t^\prime)^\ast K$ is relatively perfect over $Y_1 \times_S T$.
	Moreover, as $T$ is separated (\Cref{rem:fm_cube_separated}), \Cref{lem:left_adjoint} applied to the families $Y_i\times_S T\to T$ shows that the kernel of the integral transform which is left adjoint to $\Phi_{\mathbf{L} (t^\prime)^\ast K}$ is given by $\operatorname{\mathbf{R}\mathcal{H}\! \mathit{om}} ( \mathbf{L} (t^\prime)^\ast K, (g^\prime_2)^! \mathcal{O}_{Y_1\times_S T})$.
	It suffices to show there is an isomorphism
	\begin{displaymath}
		\operatorname{\mathbf{R}\mathcal{H}\! \mathit{om}} ( \mathbf{L} (t^\prime)^\ast K, (g^\prime_2)^! \mathcal{O}_{Y_1\times_S T})
		\to  \mathbf{L} (t^\prime)^\ast \operatorname{\mathbf{R}\mathcal{H}\! \mathit{om}} ( K, (f^\prime_2)^! \mathcal{O}_{Y_1}).
	\end{displaymath}
	By \cite[\href{https://stacks.math.columbia.edu/tag/0AU1}{Tag 0AU1}]{StacksProject}, $(f^\prime_2)^! \mathcal{O}_{Y_1} \in D^+_{\operatorname{coh}}(Y_1\times_S Y_2)$ because upper shriek functors for proper morphisms preserve complexes with bounded below and coherent cohomology.
	Then \Cref{lem:base_change_relative_dual} yields
	\begin{displaymath}
		\mathbf{L}(t^\prime)^\ast
		\operatorname{\mathbf{R}\mathcal{H}\! \mathit{om}}( K, (f^\prime_2)^! \mathcal{O}_{Y_1})
		\cong \operatorname{\mathbf{R}\mathcal{H}\! \mathit{om}}(\mathbf{L}(t^\prime)^\ast K, (g^\prime_2)^! \mathcal{O}_{Y_1\times_S T}).
	\end{displaymath}
	This completes the proof.
\end{proof}

\section{Full faithfulness \& equivalences}
\label{sec:fully_faithful_equivalence}

\begin{lemma}
	\label{lem:fm_cube_pushforward}
	Consider \Cref{setup:fm_cube}.
	Let $K\in D_{\operatorname{qc}}(Y_1\times_S Y_2)$.
	Then on $D_{\operatorname{qc}}$ there is a natural isomorphism
	\begin{displaymath}
		\beta^K \colon \Phi_K \circ \mathbf{R}(t_1)_\ast
		\to \mathbf{R}(t_2)_\ast \circ \Phi_{\mathbf{L}(t^\prime)^\ast K}.
	\end{displaymath}
\end{lemma}

\begin{proof}
	For each $E\in D_{\operatorname{qc}}(Y_1\times_S T)$, there is a string of natural isomorphisms:
	\begin{displaymath}
		\begin{aligned}
			\Phi_K \circ \mathbf{R} (t_1)_\ast (E)
			&= \mathbf{R} (f_1^\prime)_\ast ( \mathbf{L} (f^\prime_2)^\ast \mathbf{R} (t_1)_\ast E \otimes^{\mathbf{L}} K )
			\\&\cong \mathbf{R} (f_1^\prime)_\ast ( \mathbf{R} (t^\prime)_\ast \mathbf{L}(g^\prime_2)^\ast E \otimes^{\mathbf{L}} K )
			&& \textrm{(flat base change)}
			\\&\cong \mathbf{R} (f_1^\prime)_\ast \mathbf{R} (t^\prime)_\ast ( \mathbf{L}(g^\prime_2)^\ast E \otimes^{\mathbf{L}} \mathbf{L}(t^\prime)^\ast K )
			&& \textrm{(projection formula)}
			\\&\cong \mathbf{R} (t_2)_\ast \mathbf{R} (g^\prime_1)_\ast ( \mathbf{L}(g^\prime_2)^\ast E \otimes^{\mathbf{L}} \mathbf{L}(t^\prime)^\ast K )
			&& \textrm{(pseudofunctoriality)}
			\\&= \mathbf{R} (t_2)_\ast \circ \Phi_{\mathbf{L}(t^\prime)^\ast K} (E).
		\end{aligned}
	\end{displaymath}
	This completes the proof.
\end{proof}

\begin{lemma}
	\label{lem:integral_transform_mates}
	Consider \Cref{setup:fm_cube}.
	Let $K\in D_{\operatorname{qc}}(Y_1\times_S Y_2)$.
	The isomorphisms $\alpha^K$ and $\beta^K$ of \Cref{lem:fm_cube_flat_base_change,lem:fm_cube_pushforward} are respectively left and right mates under the diagram
	\begin{displaymath}
		\begin{tikzcd}
			{D_{\operatorname{qc}}(Y_1)} & {D_{\operatorname{qc}}(Y_2)} \\
			{D_{\operatorname{qc}}(Y_1\times_ST)} & {D_{\operatorname{qc}}(Y_2\times_ST).}
			\arrow["{\Phi_K}", from=1-1, to=1-2]
			\arrow[""{name=0, anchor=center, inner sep=0}, "{\mathbf{L}t_1^\ast }"', shift right=2, from=1-1, to=2-1]
			\arrow[""{name=1, anchor=center, inner sep=0}, "{\mathbf{L}t_2^\ast }"', shift right=2, from=1-2, to=2-2]
			\arrow[""{name=2, anchor=center, inner sep=0}, "{\mathbf{R}(t_1)_\ast }"', shift right=2, from=2-1, to=1-1]
			\arrow["{\Phi_{\mathbf{L}(t^\prime)^\ast K}}"', from=2-1, to=2-2]
			\arrow[""{name=3, anchor=center, inner sep=0}, "{\mathbf{R}(t_2)_\ast }"', shift right=2, from=2-2, to=1-2]
			\arrow["\dashv"{anchor=center}, draw=none, from=0, to=2]
			\arrow["\dashv"{anchor=center}, draw=none, from=1, to=3]
		\end{tikzcd}
	\end{displaymath}
\end{lemma}

\begin{remark}
	\label{rem:hall_priver}
	The compatibility of \Cref{lem:integral_transform_mates} also appears in \cite[Lemma 5.6(3)]{Hall/Priver:2024}.
	We provide a detailed proof here.
\end{remark}

\begin{proof}
	[Proof of \Cref{lem:integral_transform_mates}]
	Denote by $\xi^i\colon \mathbf{L}t_i^\ast \mathbf{R}(t_i)_\ast \to 1$ the counit of the adjunction $\mathbf{L}t_i^\ast \dashv\mathbf{R}(t_i)_\ast $ for $i=1,2$.
	By \Cref{lem:mates_equivalent_conditions}, it suffices to verify that the following diagram commutes:
	\begin{displaymath}
		\begin{tikzcd}
			{\mathbf{L}t_2^\ast\Phi_K\mathbf{R}(t_1)_\ast} & {\Phi_{\mathbf{L}(t^\prime)^\ast K}\mathbf{L}t_1^\ast \mathbf{R}(t_1)_\ast} \\
			{\mathbf{L}t_2^\ast\mathbf{R}(t_2)_\ast\Phi_{\mathbf{L}(t^\prime)^\ast K}} & {\Phi_{\mathbf{L}(t^\prime)^\ast K}.}
			\arrow["{\alpha^K_{\mathbf{R}(t_1)_\ast}}", from=1-1, to=1-2]
			\arrow["{\mathbf{L}t_2^\ast(\beta^K)}"', from=1-1, to=2-1]
			\arrow["{\Phi_{\mathbf{L}(t^\prime)^\ast K}(\xi^1)}", from=1-2, to=2-2]
			\arrow["{\xi^2_{\Phi_{\mathbf{L}(t^\prime)^\ast K}}}"', from=2-1, to=2-2]
		\end{tikzcd}
	\end{displaymath}
	For brevity, in the remainder of the proof, we drop the $\mathbf{L}$'s and $\mathbf{R}$'s in the notation for the derived functors (i.e.\ all functors are understood to be derived).
	We use the following abbreviations for the natural isomorphisms involved:
	\begin{itemize}
		\item $\operatorname{BC}_1$ and $\operatorname{BC}_2$ for the base change isomorphisms $(f_2^\prime)^\ast (t_1)_\ast \to (t^\prime)_\ast (g_2^\prime)^\ast$ and $t_2^\ast (f_1^\prime)_\ast\to (g_1^\prime)_\ast (t^\prime)^\ast$, respectively
		\item $\operatorname{PF}$ for the projection formula
		\item $F_{(-)_\ast}$ and $F_{(-)^\ast}$ for the pseudofunctoriality isomorphisms of $(-)_\ast$ and of $(-)^\ast$, respectively
		\item $M_{(-)^\ast}$ for monoidality of $(-)^\ast$
		\item $\xi^\prime$ for the counit of $(t^\prime)^\ast \dashv (t^\prime)_\ast$.
	\end{itemize}
	Choose $E\in D_{\operatorname{qc}}(Y_1\times_S T)$.
	By the proofs of \Cref{lem:fm_cube_flat_base_change,lem:fm_cube_pushforward}, the composite along the top and right of the square above, evaluated at $E$, is the composite of $\operatorname{BC}_2$, $M_{(-)^\ast}$, $F_{(-)^\ast}$, and $\xi^1$, whereas the composite along the left and bottom is the composite of $\operatorname{BC}_1$, $\operatorname{PF}$, $F_{(-)_\ast}$, and $\xi^2$.
	Consider the following diagrams:
	\begin{equation}
		\label{diag:pastingA}
		\begin{tikzcd}
			{t_2^\ast (f_1^\prime)_\ast ((f_2^\prime)^\ast (t_1)_\ast E\otimes K)} & {t_2^\ast (f_1^\prime)_\ast ((t^\prime)_\ast (g_2^\prime)^\ast E\otimes K)} \\
			{(g_1^\prime)_\ast (t^\prime)^\ast ((f_2^\prime)^\ast (t_1)_\ast E\otimes K)} & {(g_1^\prime)_\ast (t^\prime)^\ast ((t^\prime)_\ast (g_2^\prime)^\ast E\otimes K)}
			\arrow["{\operatorname{BC}_1}", from=1-1, to=1-2]
			\arrow["{\operatorname{BC}_2}"', from=1-1, to=2-1]
			\arrow["{\operatorname{BC}_2}", from=1-2, to=2-2]
			\arrow["{\operatorname{BC}_1}"', from=2-1, to=2-2]
		\end{tikzcd}
	\end{equation}
	\begin{equation}
		\label{diag:pastingB}
		\begin{tikzcd}
			{t_2^\ast (f_1^\prime)_\ast ((t^\prime)_\ast (g_2^\prime)^\ast E\otimes K)} & {t_2^\ast (f_1^\prime)_\ast (t^\prime)_\ast ((g_2^\prime)^\ast E\otimes (t^\prime)^\ast K)} \\
			{(g_1^\prime)_\ast (t^\prime)^\ast ((t^\prime)_\ast (g_2^\prime)^\ast E\otimes K)} & {(g_1^\prime)_\ast (t^\prime)^\ast (t^\prime)_\ast ((g_2^\prime)^\ast E\otimes (t^\prime)^\ast K)}
			\arrow["{\operatorname{PF}}", from=1-1, to=1-2]
			\arrow["{\operatorname{BC}_2}"', from=1-1, to=2-1]
			\arrow["{\operatorname{BC}_2}", from=1-2, to=2-2]
			\arrow["{\operatorname{PF}}"', from=2-1, to=2-2]
		\end{tikzcd}
	\end{equation}
	\begin{equation}
		\label{diag:pastingC}
		\begin{tikzcd}
			{(g_1^\prime)_\ast (t^\prime)^\ast ((f_2^\prime)^\ast (t_1)_\ast E\otimes K)} & {(g_1^\prime)_\ast (t^\prime)^\ast ((t^\prime)_\ast (g_2^\prime)^\ast E\otimes K)} \\
			{(g_1^\prime)_\ast ((t^\prime)^\ast (f_2^\prime)^\ast (t_1)_\ast E\otimes (t^\prime)^\ast K)} & {(g_1^\prime)_\ast ((t^\prime)^\ast (t^\prime)_\ast (g_2^\prime)^\ast E\otimes (t^\prime)^\ast K)}
			\arrow["{\operatorname{BC}_1}", from=1-1, to=1-2]
			\arrow["{M_{(-)^\ast}}"', from=1-1, to=2-1]
			\arrow["{M_{(-)^\ast}}", from=1-2, to=2-2]
			\arrow["{\operatorname{BC}_1}"', from=2-1, to=2-2]
		\end{tikzcd}
	\end{equation}
	\begin{equation}
		\label{diag:pastingD}
		\begin{tikzcd}
			{(g_1^\prime)_\ast (t^\prime)^\ast ((t^\prime)_\ast (g_2^\prime)^\ast E\otimes K)} & {(g_1^\prime)_\ast (t^\prime)^\ast (t^\prime)_\ast ((g_2^\prime)^\ast E\otimes (t^\prime)^\ast K)} \\
			{(g_1^\prime)_\ast ((t^\prime)^\ast (t^\prime)_\ast (g_2^\prime)^\ast E\otimes (t^\prime)^\ast K)} & {(g_1^\prime)_\ast ((g_2^\prime)^\ast E\otimes (t^\prime)^\ast K)}
			\arrow["{\operatorname{PF}}", from=1-1, to=1-2]
			\arrow["{M_{(-)^\ast}}"', from=1-1, to=2-1]
			\arrow["{\xi^\prime}", from=1-2, to=2-2]
			\arrow["{\xi^\prime}"', from=2-1, to=2-2]
		\end{tikzcd}
	\end{equation}
	\begin{equation}
		\label{diag:pastingE}
		\begin{tikzcd}
			{t_2^\ast (f_1^\prime)_\ast (t^\prime)_\ast ((g_2^\prime)^\ast E\otimes (t^\prime)^\ast K)} & {t_2^\ast (t_2)_\ast (g_1^\prime)_\ast ((g_2^\prime)^\ast E\otimes (t^\prime)^\ast K)} \\
			{(g_1^\prime)_\ast (t^\prime)^\ast (t^\prime)_\ast ((g_2^\prime)^\ast E\otimes (t^\prime)^\ast K)} & {(g_1^\prime)_\ast ((g_2^\prime)^\ast E\otimes (t^\prime)^\ast K)}
			\arrow["{F_{(-)_\ast}}", from=1-1, to=1-2]
			\arrow["{\operatorname{BC}_2}"', from=1-1, to=2-1]
			\arrow["{\xi^2}", from=1-2, to=2-2]
			\arrow["{\xi^\prime}"', from=2-1, to=2-2]
		\end{tikzcd}
	\end{equation}
	\begin{equation}
		\label{diag:pastingF}
		\begin{tikzcd}
			{(g_1^\prime)_\ast ((t^\prime)^\ast (f_2^\prime)^\ast (t_1)_\ast E\otimes (t^\prime)^\ast K)} & {(g_1^\prime)_\ast ((t^\prime)^\ast (t^\prime)_\ast (g_2^\prime)^\ast E\otimes (t^\prime)^\ast K)} \\
			{(g_1^\prime)_\ast ((g_2^\prime)^\ast t_1^\ast (t_1)_\ast E\otimes (t^\prime)^\ast K)} & {(g_1^\prime)_\ast ((g_2^\prime)^\ast E\otimes (t^\prime)^\ast K)}
			\arrow["{\operatorname{BC}_1}", from=1-1, to=1-2]
			\arrow["{F_{(-)^\ast}}"', from=1-1, to=2-1]
			\arrow["{\xi^\prime}", from=1-2, to=2-2]
			\arrow["{\xi^1}"', from=2-1, to=2-2]
		\end{tikzcd}
	\end{equation}
	Assume that \eqref{diag:pastingA}, \eqref{diag:pastingB}, \eqref{diag:pastingC}, \eqref{diag:pastingD}, \eqref{diag:pastingE}, and \eqref{diag:pastingF} commute.
	Then, using \eqref{diag:pastingF}, \eqref{diag:pastingC}, \eqref{diag:pastingA}, \eqref{diag:pastingD}, \eqref{diag:pastingB}, and \eqref{diag:pastingE} in this order, the composite of $\operatorname{BC}_2$, $M_{(-)^\ast}$, $F_{(-)^\ast}$, and $\xi^1$ is successively equal to the composites of
	\begin{displaymath}
		\begin{aligned}
			&\operatorname{BC}_2,\ M_{(-)^\ast},\ \operatorname{BC}_1,\ \xi^\prime,
			&&\operatorname{BC}_2,\ \operatorname{BC}_1,\ M_{(-)^\ast},\ \xi^\prime,
			&&\operatorname{BC}_1,\ \operatorname{BC}_2,\ M_{(-)^\ast},\ \xi^\prime,
			\\&\operatorname{BC}_1,\ \operatorname{BC}_2,\ \operatorname{PF},\ \xi^\prime,
			&&\operatorname{BC}_1,\ \operatorname{PF},\ \operatorname{BC}_2,\ \xi^\prime,
			&&\operatorname{BC}_1,\ \operatorname{PF},\ F_{(-)_\ast},\ \xi^2,
		\end{aligned}
	\end{displaymath}
	which gives the desired claim.

	We explain why the diagrams commute.
	For \eqref{diag:pastingA} and \eqref{diag:pastingB}, use the naturality of $\operatorname{BC}_2$, whereas \eqref{diag:pastingC} commutes by the naturality of $M_{(-)^\ast}$.
	The diagram \eqref{diag:pastingD} is the functor $(g_1^\prime)_\ast$ applied to the following diagram evaluated at $((g_2^\prime)^\ast E,K)$:
	\begin{displaymath}
		\begin{tikzcd}
			{(t^\prime)^\ast ((t^\prime)_\ast (-)\otimes -)} & {(t^\prime)^\ast (t^\prime)_\ast (-\otimes (t^\prime)^\ast (-))} \\
			{(t^\prime)^\ast (t^\prime)_\ast (-)\otimes (t^\prime)^\ast(-)} & {(-)\otimes (t^\prime)^\ast(-)}
			\arrow["{\operatorname{PF}}", from=1-1, to=1-2]
			\arrow["{M_{(-)^\ast}}"', from=1-1, to=2-1]
			\arrow["{\xi^\prime}", from=1-2, to=2-2]
			\arrow["{\xi^\prime}"', from=2-1, to=2-2]
		\end{tikzcd}
	\end{displaymath}
	which commutes by definition of the projection formula, see e.g.\ the proof of \cite[Proposition 22.81]{Gortz/Wedhorn:2023}.
	The diagram \eqref{diag:pastingE} is the following diagram evaluated at $(g_2^\prime)^\ast E\otimes (t^\prime)^\ast K$:
	\begin{displaymath}
		\begin{tikzcd}
			{t_2^\ast (f^\prime_1)_\ast (t^\prime)_\ast} & {t_2^\ast (t_2)_\ast (g^\prime_1)_\ast} \\
			{(g^\prime_1)_\ast (t^\prime)^\ast (t^\prime)_\ast} & {(g_1^\prime)_\ast}
			\arrow["{F_{(-)_\ast}}", from=1-1, to=1-2]
			\arrow["{\operatorname{BC}_2}"', from=1-1, to=2-1]
			\arrow["{\xi^2}", from=1-2, to=2-2]
			\arrow["{\xi^\prime}"', from=2-1, to=2-2]
		\end{tikzcd}
	\end{displaymath}
	which commutes by definition of the base change morphism, see \cite[Definition and Remark 21.129]{Gortz/Wedhorn:2023} or \cite[(3.7.2)]{Lipman:2009}.
	Lastly, \eqref{diag:pastingF} is the functor $(g_1^\prime)_\ast (-\otimes (t^\prime)^\ast K)$ applied to the following diagram evaluated at $E$:
	\begin{displaymath}
		\begin{tikzcd}
			{(t^\prime)^\ast (f_2^\prime)^\ast (t_1)_\ast} & {(t^\prime)^\ast (t^\prime)_\ast (g_2^\prime)^\ast} \\
			{(g_2^\prime)^\ast t_1^\ast (t_1)_\ast} & {(g_2^\prime)^\ast}
			\arrow["{\operatorname{BC}_1}", from=1-1, to=1-2]
			\arrow["{F_{(-)^\ast}}"', from=1-1, to=2-1]
			\arrow["{\xi^\prime}", from=1-2, to=2-2]
			\arrow["{\xi^1}"', from=2-1, to=2-2]
		\end{tikzcd}
	\end{displaymath}
	which commutes by definition of the base change morphism.
\end{proof}

\begin{remark}
	\label{rem:units_compatibility}
	Assume the hypothesis of \Cref{lem:integral_transform_mates}.
	Denote by $\zeta^i$ the unit of $\mathbf{L}t_i^\ast \dashv \mathbf{R}(t_i)_\ast$ for $i\in\{1,2\}$.
	From \Cref{lem:mates_equivalent_conditions,lem:integral_transform_mates}, the following diagram commutes:
	\begin{displaymath}
		\begin{tikzcd}
			{\Phi_K} & {\mathbf{R}(t_2)_\ast \mathbf{L}t_2^\ast \Phi_K} \\
			{\Phi_K \mathbf{R}(t_1)_\ast \mathbf{L}t_1^\ast} & {\mathbf{R}(t_2)_\ast \Phi_{\mathbf{L}(t^\prime)^\ast K}\mathbf{L}t_1^\ast.}
			\arrow["{\zeta^2_{\Phi_K}}", from=1-1, to=1-2]
			\arrow["{\Phi_K(\zeta^1)}"', from=1-1, to=2-1]
			\arrow["{\mathbf{R}(t_2)_\ast (\alpha^K)}", from=1-2, to=2-2]
			\arrow["{\beta^K_{\mathbf{L}t_1^\ast}}"', from=2-1, to=2-2]
		\end{tikzcd}
	\end{displaymath}
\end{remark}

The following does not follow formally from \Cref{lem:fm_cube_flat_base_change}, as the mate of an isomorphism need not be an isomorphism, see \Cref{rem:mates_of_isomorphisms}.

\begin{proposition}
	\label{prop:left_mate_isomorphism}
	Consider \Cref{setup:fm_cube} where $t$ is affine.
	Let $K\in D^b_{\operatorname{coh}}(Y_1\times_S Y_2)$ be relatively perfect over $Y_2$.
	Denote by $K^\prime$ the kernel obtained in \Cref{prop:right_adjoint_base_change}.
	Consider the diagram:
	\begin{displaymath}
		\begin{tikzcd}[column sep=large]
			{D_{\operatorname{qc}}(Y_1)} & {D_{\operatorname{qc}}(Y_1\times_S T)} \\
			{D_{\operatorname{qc}}(Y_2)} & {D_{\operatorname{qc}}(Y_2\times_S T).}
			\arrow["{\mathbf{L}t_1^\ast}", from=1-1, to=1-2]
			\arrow[""{name=0, anchor=center, inner sep=0}, "{\Phi_K}"', shift right=2, from=1-1, to=2-1]
			\arrow[""{name=1, anchor=center, inner sep=0}, "{\Phi_{\mathbf{L}(t^\prime)^\ast K}}"', shift right=2, from=1-2, to=2-2]
			\arrow[""{name=2, anchor=center, inner sep=0}, "{\Phi_{K^\prime}}"', shift right=2, from=2-1, to=1-1]
			\arrow["{\mathbf{L}t_2^\ast}"', from=2-1, to=2-2]
			\arrow[""{name=3, anchor=center, inner sep=0}, "{\Phi_{\mathbf{L}(t^\prime)^\ast K^\prime}}"', shift right=2, from=2-2, to=1-2]
			\arrow["\dashv"{anchor=center}, draw=none, from=0, to=2]
			\arrow["\dashv"{anchor=center}, draw=none, from=1, to=3]
		\end{tikzcd}
	\end{displaymath}
	Then, under this diagram, the isomorphism of \Cref{lem:fm_cube_flat_base_change} (applied with the roles of $Y_1$ and $Y_2$ interchanged)
	\begin{displaymath}
		\alpha^{K^\prime}\colon \mathbf{L}t_1^\ast\Phi_{K^\prime}
		\to \Phi_{\mathbf{L}(t^\prime)^\ast K^\prime} \mathbf{L}t_2^\ast
	\end{displaymath}
	has a left mate
	\begin{displaymath}
		\lambda^K\colon \Phi_{\mathbf{L}(t^\prime)^\ast K} \mathbf{L}t_1^\ast \to \mathbf{L}t_2^\ast\Phi_K
	\end{displaymath}
	that is an isomorphism.
\end{proposition}

\begin{proof}
	Denote by $\eta$ and $\epsilon$ the unit and counit of $\Phi_K\dashv\Phi_{K^\prime}$.
	Denote by $\eta^\prime$ and $\epsilon^\prime$ the unit and counit of $\Phi_{\mathbf{L}(t^\prime)^\ast K} \dashv \Phi_{\mathbf{L}(t^\prime)^\ast K^\prime}$.
	By \Cref{prop:right_adjoint_base_change}, the latter adjunction exists.
	By \Cref{def:mates}, the component of $\lambda^K$ at $E\in D_{\operatorname{qc}}(Y_1)$ is the composite
	\begin{displaymath}
		\begin{aligned}
			\Phi_{\mathbf{L}(t^\prime)^\ast K}
			\mathbf{L}t_1^\ast E
			&\xrightarrow{\Phi_{\mathbf{L}(t^\prime)^\ast K} \mathbf{L}t_1^\ast(\eta_E)}
			\Phi_{\mathbf{L}(t^\prime)^\ast K}
			\mathbf{L}t_1^\ast\Phi_{K^\prime}\Phi_K(E)
			\\&\xrightarrow{\Phi_{\mathbf{L}(t^\prime)^\ast K}(\alpha^{K^\prime}_{\Phi_K(E)})}
			\Phi_{\mathbf{L}(t^\prime)^\ast K}
			\Phi_{\mathbf{L}(t^\prime)^\ast K^\prime}
			\mathbf{L}t_2^\ast\Phi_K(E)
			\\&\xrightarrow{\epsilon^\prime_{\mathbf{L}t_2^\ast\Phi_K(E)}}
			\mathbf{L}t_2^\ast\Phi_K(E).
		\end{aligned}
	\end{displaymath}
	Choose $E\in D_{\operatorname{qc}}(Y_1)$ and $A\in D_{\operatorname{qc}}(Y_2\times_S T)$.
	Precomposition with $\lambda^K_E$ gives a morphism
	\begin{displaymath}
		\operatorname{Hom} (\mathbf{L}t_2^\ast\Phi_K(E),A)
		\to \operatorname{Hom} ( \Phi_{\mathbf{L}(t^\prime)^\ast K} \mathbf{L}t_1^\ast E, A ).
	\end{displaymath}
	We claim that this morphism is the composite
	\begin{equation}
		\label{eq:adjunctions}
		\begin{aligned}
			\operatorname{Hom}(\mathbf{L}t_2^\ast\Phi_K(E),A)
			&\cong \operatorname{Hom} (\Phi_K(E),\mathbf{R}(t_2)_\ast A)
			\\&\cong \operatorname{Hom} (E,\Phi_{K^\prime}\mathbf{R}(t_2)_\ast A)
			\\&\xrightarrow{\beta^{K^\prime}_A\circ-}\operatorname{Hom}(E, \mathbf{R}(t_1)_\ast \Phi_{\mathbf{L}(t^\prime)^\ast K^\prime}(A) )
			\\&\cong \operatorname{Hom} (\mathbf{L}t_1^\ast E, \Phi_{\mathbf{L}(t^\prime)^\ast K^\prime}(A))
			\\&\cong \operatorname{Hom} ( \Phi_{\mathbf{L}(t^\prime)^\ast K} \mathbf{L}t_1^\ast E, A),
		\end{aligned}
	\end{equation}
	where $\beta^{K^\prime}$ is the natural isomorphism of \Cref{lem:fm_cube_pushforward} applied to $K^\prime$ (with the roles of $Y_1$ and $Y_2$ interchanged).

	We verify the claim.
	Choose a morphism $h\colon \mathbf{L}t_2^\ast\Phi_K(E)\to A$, and let $h^\flat \colon \Phi_K(E)\to\mathbf{R}(t_2)_\ast A$ be its adjunct with respect to $\mathbf{L}t_2^\ast\dashv\mathbf{R}(t_2)_\ast$.
	Under the adjunction $\Phi_{\mathbf{L}(t^\prime)^\ast K} \dashv \Phi_{\mathbf{L}(t^\prime)^\ast K^\prime}$, the morphism $h\circ\lambda^K_E$ corresponds to
	\begin{equation}
		\label{eq:adjunctions2}
		\mathbf{L}t_1^\ast E
		\xrightarrow{\mathbf{L}t_1^\ast(\eta_E)}
		\mathbf{L}t_1^\ast\Phi_{K^\prime}\Phi_K(E)
		\xrightarrow{\alpha^{K^\prime}_{\Phi_K(E)}}
		\Phi_{\mathbf{L}(t^\prime)^\ast K^\prime}
		\mathbf{L}t_2^\ast\Phi_K(E)
		\xrightarrow{ \Phi_{\mathbf{L}(t^\prime)^\ast K^\prime}(h) }
		\Phi_{\mathbf{L}(t^\prime)^\ast K^\prime}(A).
	\end{equation}
	By \Cref{lem:integral_transform_mates} applied to the kernel $K^\prime$ (with the roles of $Y_1$ and $Y_2$ interchanged), $\alpha^{K^\prime}$ and $\beta^{K^\prime}$ are mates with respect to $\mathbf{L}t_2^\ast\dashv\mathbf{R}(t_2)_\ast$ and $\mathbf{L}t_1^\ast\dashv\mathbf{R}(t_1)_\ast$.
	Then \Cref{lem:mates_equivalent_conditions}\eqref{lem:mates_equivalent_conditions4}, applied with $c=\Phi_K(E)$ and $d=A$, says that the adjunct of $\Phi_{\mathbf{L}(t^\prime)^\ast K^\prime}(h) \circ \alpha^{K^\prime}_{\Phi_K(E)}$ is
	\begin{displaymath}
		\beta^{K^\prime}_A
		\circ \Phi_{K^\prime}(h^\flat) \colon
		\Phi_{K^\prime}\Phi_K(E)
		\to \mathbf{R}(t_1)_\ast \Phi_{\mathbf{L}(t^\prime)^\ast K^\prime}(A).
	\end{displaymath}
	After precomposing with the unit $\eta_E \colon E\to\Phi_{K^\prime}\Phi_K(E)$, the adjunct of \eqref{eq:adjunctions2} is
	\begin{displaymath}
		E
		\xrightarrow{\eta_E} \Phi_{K^\prime}\Phi_K(E)
		\xrightarrow{\Phi_{K^\prime}(h^\flat)} \Phi_{K^\prime}\mathbf{R}(t_2)_\ast A
		\xrightarrow{\beta^{K^\prime}_A} \mathbf{R}(t_1)_\ast \Phi_{\mathbf{L}(t^\prime)^\ast K^\prime}(A).
	\end{displaymath}
	This is the image of $h$ under the composite \eqref{eq:adjunctions}.
	Hence, precomposition with $\lambda^K_E$ agrees with \eqref{eq:adjunctions}.

	Every arrow in \eqref{eq:adjunctions} is a bijection.
	Indeed, the first, second, fourth, and fifth arrows are adjunction isomorphisms, while the third is induced by the natural isomorphism $\beta^{K^\prime}$.
	Consequently, precomposition with $\lambda^K_E$ induces a bijection
	\begin{displaymath}
		\operatorname{Hom} (\mathbf{L}t_2^\ast\Phi_K(E),A)
		\to
		\operatorname{Hom}( \Phi_{\mathbf{L}(t^\prime)^\ast K} \mathbf{L}t_1^\ast E, A)
	\end{displaymath}
	for every $A\in D_{\operatorname{qc}}(Y_2\times_S T)$.
	By \cite[\href{https://stacks.math.columbia.edu/tag/001P}{Tag 001P}]{StacksProject}, $\lambda^K_E$ is an isomorphism.
	Since this holds for every $E\in D_{\operatorname{qc}}(Y_1)$, $\lambda^K$ is a natural isomorphism.
\end{proof}

\begin{remark}
	\label{rem:unit_counit_compatibility}
	Consider the situation of \Cref{prop:left_mate_isomorphism}.
	By \Cref{lem:mates_equivalent_conditions}, it follows that
	\begin{equation}
		\label{eq:unit_compatibility}
		\Phi_{\mathbf{L}(t^\prime)^\ast K^\prime} (\lambda^K) \circ \eta^\prime \mathbf{L}t^\ast_1 = \alpha^{K^\prime} \Phi_K \circ \mathbf{L}t^\ast_1 (\eta)
	\end{equation}
	and
	\begin{equation}
		\label{eq:counit_compatibility}
		\mathbf{L}t^\ast_2 (\epsilon) \circ \lambda^K\Phi_{K^\prime} = \epsilon^\prime \mathbf{L}t_2^\ast \circ \Phi_{\mathbf{L}(t^\prime)^\ast K} (\alpha^{K^\prime}).
	\end{equation}
	As $\lambda^K$ and $\alpha^{K^\prime}$ are isomorphisms, \eqref{eq:unit_compatibility} identifies the derived pullback $\mathbf{L}t_1^\ast(\eta)$ of the unit of $\Phi_K\dashv \Phi_{K^\prime}$ with the unit $\eta^\prime\mathbf{L}t_1^\ast$ of the base changed adjunction, and \eqref{eq:counit_compatibility} does the same for the counits.
\end{remark}

\begin{lemma}
	\label{lem:equivalences_induced}
	Let $f_1\colon Y_1\to S$ and $f_2\colon Y_2\to S$ be proper flat morphisms to a separated Noetherian scheme.
	Suppose $K\in D^b_{\operatorname{coh}}(Y_1\times_S Y_2)$ is relatively perfect over each $Y_i$.
	Then the following are equivalent:
	\begin{enumerate}
		\item \label{lem:equivalences_induced1} $\Phi_K$ is fully faithful (resp.\ an equivalence) on $D_{\operatorname{qc}}$
		\item \label{lem:equivalences_induced2} $\Phi_K$ restricts to a fully faithful functor (resp.\ an equivalence) on $D^b_{\operatorname{coh}}$
		\item \label{lem:equivalences_induced3} $\Phi_K$ restricts to a fully faithful functor (resp.\ an equivalence) on $\operatorname{Perf}$.
	\end{enumerate}
	In particular, full faithfulness can be tested on a compact generator of $D_{\operatorname{qc}}(Y_1)$, while equivalences can be tested on compact generators for each $D_{\operatorname{qc}}(Y_i)$.
\end{lemma}

\begin{proof}
	Denote by $K^\prime$ the kernel of the integral transform obtained in \Cref{prop:right_adjoint_base_change} (with $t=1_S$) which is right adjoint to $\Phi_K$ on $D_{\operatorname{qc}}$.
	Denote by $K^{\prime \prime}$ the kernel of the integral transform obtained in \Cref{prop:left_adjoint_base_change} (with $t=1_S$) which is left adjoint to $\Phi_K$ on $D_{\operatorname{qc}}$.
	We start with full faithfulness.
	Since $K$ is relatively perfect over each $Y_i$, $\Phi_K$ restricts to functors on $D^b_{\operatorname{coh}}$ and $\operatorname{Perf}$ \cite[Propositions 3.2 and 3.6]{Dutta/Lank/ManaliRahul:2025}.
	Hence, $\eqref{lem:equivalences_induced1} \implies \eqref{lem:equivalences_induced2} \implies \eqref{lem:equivalences_induced3}$.

	Assume \eqref{lem:equivalences_induced3} holds; that is, $\Phi_K$ restricts to a fully faithful functor on $\operatorname{Perf}$.
	Denote by $e\colon \Phi_{K^{\prime \prime}} \circ \Phi_K \to 1$ the counit of the adjoint pair $\Phi_{K^{\prime \prime}}$ and $\Phi_K$.
	As $Y_1$ is separated (\Cref{rem:fm_cube_separated}), \Cref{lem:relative_dual_relatively_perfect} implies $K^{\prime \prime}$ is relatively perfect over $Y_1$.
	As $Y_1$ is the target of $\Phi_{K^{\prime\prime}}$, \cite[Propositions 3.2 and 3.6]{Dutta/Lank/ManaliRahul:2025} imply $\Phi_{K^{\prime \prime}}$ restricts to a functor on $\operatorname{Perf}$.
	Since $\Phi_K$ restricts to a fully faithful functor on $\operatorname{Perf}$, \cite[\href{https://stacks.math.columbia.edu/tag/07RB}{Tag 07RB}]{StacksProject} implies $e_P$ is an isomorphism for all $P\in \operatorname{Perf}(Y_1)$.
	Consider the strictly full subcategory $\mathcal{T}$ of objects $B\in D_{\operatorname{qc}}(Y_1)$ such that $e_B$ is an isomorphism.
	Since both integral transforms preserve small coproducts, $\mathcal{T}$ is a localizing subcategory containing $\operatorname{Perf}(Y_1)$.
	As $D_{\operatorname{qc}}(Y_1)$ is compactly generated, \Cref{lem:localizing_iff_cpt_gen} implies $D_{\operatorname{qc}}(Y_1) = \mathcal{T}$.
	Then \cite[\href{https://stacks.math.columbia.edu/tag/07RB}{Tag 07RB}]{StacksProject} implies $\Phi_K$ is fully faithful on $D_{\operatorname{qc}}$.

	Next we prove the case for equivalences.
	By \Cref{lem:equivalence_via_compacts}, $\eqref{lem:equivalences_induced1} \iff \eqref{lem:equivalences_induced3}$.
	By \Cref{cor:right_adjoint_restricts}, $\Phi_K$ and $\Phi_{K^\prime}$ restrict to an adjoint pair on $D^b_{\operatorname{coh}}$.
	Assume \eqref{lem:equivalences_induced1} holds.
	Then the unit and counit of $\Phi_K\dashv\Phi_{K^\prime}$ are isomorphisms, in particular on $D^b_{\operatorname{coh}}$, and so \cite[\href{https://stacks.math.columbia.edu/tag/07RB}{Tag 07RB}]{StacksProject} yields \eqref{lem:equivalences_induced2}.
	Assume \eqref{lem:equivalences_induced2} holds.
	Since $\Phi_K$ restricts to a functor on $\operatorname{Perf}$, it is a fully faithful functor on $\operatorname{Perf}$.
	By the fully faithful case, $\Phi_K$ is fully faithful on $D_{\operatorname{qc}}$.
	Let $\mathcal{S}$ be the strictly full subcategory of $B\in D_{\operatorname{qc}}(Y_2)$ such that the counit $\epsilon_B \colon (\Phi_K \circ \Phi_{K^\prime})(B) \to B$ is an isomorphism.
	Then $\mathcal{S}$ is localizing.
	Since $\epsilon_B$ is a morphism of $D^b_{\operatorname{coh}}(Y_2)$ for all $B\in D^b_{\operatorname{coh}}(Y_2)$, the assumption implies $D^b_{\operatorname{coh}}(Y_2) \subseteq \mathcal{S}$.
	Thus, $\operatorname{Perf}(Y_2)\subseteq \mathcal{S}$, and so \Cref{lem:localizing_iff_cpt_gen} implies $\mathcal{S}=D_{\operatorname{qc}}(Y_2)$.
	This shows that $\eqref{lem:equivalences_induced2} \implies \eqref{lem:equivalences_induced1}$.

	We prove the last claim.
	Denote by $\eta$ and $\epsilon$ the unit and counit of $\Phi_K\dashv\Phi_{K^\prime}$.
	For $i=1,2$, let $G_i$ be a compact generator of $D_{\operatorname{qc}}(Y_i)$, see \Cref{ex:perfect_compacts_classical_generator}.
	Then $\Phi_K$ is fully faithful (resp.\ an equivalence) if, and only if, $\eta_{G_1}$ is an isomorphism (resp.\ both $\eta_{G_1}$ and $\epsilon_{G_2}$ are isomorphisms).
	Indeed, the strictly full subcategories of $D_{\operatorname{qc}}(Y_1)$ and $D_{\operatorname{qc}}(Y_2)$ consisting of objects where $\eta$ and $\epsilon$, respectively, are isomorphisms are closed under shifts, iterated cones, and small coproducts (hence, homotopy colimits).
	Therefore, \cite[\href{https://stacks.math.columbia.edu/tag/09SN}{Tag 09SN}]{StacksProject} tells us that if $\eta_{G_1}$ (resp.\ $\epsilon_{G_2}$) is an isomorphism, then $\eta$ (resp.\ $\epsilon$) is an isomorphism.
	This completes the proof.
\end{proof}

\begin{lemma}
	\label{lem:fiberwise_vanishing}
	Let $f\colon Y\to S$ be a morphism of Noetherian schemes which maps closed points to closed points.
	For a closed point $s\in S$, denote by $j_s\colon Y\times_S\operatorname{Spec}(\kappa(s))\to Y$ the canonical morphism.
	If $E\in D_{\operatorname{coh}}^b(Y)$ satisfies $\mathbf{L}j_s^\ast E\cong 0$ for every closed point $s\in S$, then $E\cong 0$.
\end{lemma}

\begin{proof}
	It suffices to show that $\operatorname{Supp}(E)$ is empty.
	Since $\operatorname{Supp}(E)$ is closed, if it were nonempty, then it would contain a closed point of $Y$ \cite[Exercise 3.13]{Gortz/Wedhorn:2020}.
	Hence, it suffices to show $\operatorname{Supp}(E)$ contains no closed points.
	Let $q\in Y$ be a closed point.
	Then $s:=f(q)$ is a closed point of $S$.
	The composite $f\circ i_q\colon \operatorname{Spec}(\kappa(q))\to S$ has image $s$, and so factors through $\operatorname{Spec}(\kappa(s))\to S$ \cite[\href{https://stacks.math.columbia.edu/tag/01J5}{Tag 01J5}]{StacksProject}.
	By the universal property of the fiber product, there is a morphism $b\colon \operatorname{Spec}(\kappa(q))\to Y\times_S\operatorname{Spec}(\kappa(s))$ such that $i_q = j_s\circ b$.
	Therefore, $\mathbf{L}i_q^\ast E\cong \mathbf{L}b^\ast\mathbf{L}j_s^\ast E\cong 0$.
	By \cite[Lemma A.1]{Iyengar/Lipman/Neeman:2015}, $q\not\in \operatorname{Supp}(E)$, as desired.
\end{proof}

\begin{proposition}
	\label{prop:descent}
	Consider \Cref{setup:fm_cube} with $S$ separated.
	Let $K\in D^b_{\operatorname{coh}}(Y_1\times_S Y_2)$ be relatively perfect over each $Y_i$.
	Suppose one of the following holds:
	\begin{enumerate}
		\item \label{prop:descent1} $t$ is affine and faithfully flat, and $\Phi_{\mathbf{L}(t^\prime)^\ast K}$ is fully faithful (resp.\ an equivalence) on $D^b_{\operatorname{coh}}$
		\item \label{prop:descent2} $\Phi_{\mathbf{L}(t^\prime)^\ast K}$ is fully faithful (resp.\ an equivalence) on $D^b_{\operatorname{coh}}$ for every closed immersion $t\colon \operatorname{Spec}(\kappa(p))\to S$ associated to a closed point $p\in S$.
	\end{enumerate}
	Then $\Phi_K$ is fully faithful (resp.\ an equivalence) on $D^b_{\operatorname{coh}}$.
\end{proposition}

\begin{proof}
	We prove the claim for full faithfulness.
	In both cases, $t$ is affine, and so $T$ is separated (\Cref{rem:fm_cube_separated}).
	Denote by $K^\prime$ the kernel obtained in \Cref{prop:right_adjoint_base_change}.
	Denote by $\eta$ (resp.\ $\eta^\prime$) the unit of the adjunction $\Phi_K\dashv\Phi_{K^\prime}$ (resp.\ $\Phi_{\mathbf{L}(t^\prime)^\ast K}\dashv\Phi_{\mathbf{L}(t^\prime)^\ast K^\prime}$) on $D_{\operatorname{qc}}$.
	By \Cref{lem:relative_perfectness_affine_base_change}, $\mathbf{L}(t^\prime)^\ast K$ is relatively perfect over each $Y_i\times_S T$.
	Hence, by \Cref{cor:right_adjoint_restricts} (applied over $S$ and, using \eqref{eq:right_adjoint_kernel_base_change}, over $T$), both adjunctions restrict to $D^b_{\operatorname{coh}}$.
	By \cite[\href{https://stacks.math.columbia.edu/tag/07RB}{Tag 07RB}]{StacksProject}, it suffices to show that $\eta_E$ is an isomorphism for all $E\in D^b_{\operatorname{coh}}(Y_1)$.
	Let $E\in D^b_{\operatorname{coh}}(Y_1)$.
	Then $\operatorname{cone}(\eta_E)\in D^b_{\operatorname{coh}}(Y_1)$.
	As $\lambda^K$ and $\alpha^{K^\prime}$ are isomorphisms (\Cref{prop:left_mate_isomorphism}), \eqref{eq:unit_compatibility} implies
	\begin{displaymath}
		\operatorname{cone}(\eta^\prime_{\mathbf{L}t_1^\ast E})\cong \mathbf{L}t_1^\ast\operatorname{cone}(\eta_E).
	\end{displaymath}

	We start with \eqref{prop:descent1}.
	Since $t$ is faithfully flat, so is $t_1$, and hence $\mathbf{L}t_1^\ast E\in D^b_{\operatorname{coh}}(Y_1\times_S T)$.
	As $\Phi_{\mathbf{L}(t^\prime)^\ast K}$ restricts to a fully faithful functor on $D^b_{\operatorname{coh}}$, $\eta^\prime_{\mathbf{L}t_1^\ast E}$ is an isomorphism \cite[\href{https://stacks.math.columbia.edu/tag/07RB}{Tag 07RB}]{StacksProject}.
	Thus, $\mathbf{L}t_1^\ast\operatorname{cone}(\eta_E)\cong 0$, and so \Cref{lem:faithfully_flat_conservative} implies $\operatorname{cone}(\eta_E)\cong 0$.

	Next, we check \eqref{prop:descent2}.
	For every closed immersion $t\colon \operatorname{Spec}(\kappa(p))\to S$ associated to a closed point $p\in S$, \Cref{lem:equivalences_induced} implies $\Phi_{\mathbf{L}(t^\prime)^\ast K}$ is fully faithful on $D_{\operatorname{qc}}$.
	Hence, $\eta^\prime$ is an isomorphism on $D_{\operatorname{qc}}$, and so $\mathbf{L}t_1^\ast\operatorname{cone}(\eta_E)\cong 0$.
	As $f_1$ is proper, it maps closed points to closed points.
	By \Cref{lem:fiberwise_vanishing}, $\operatorname{cone}(\eta_E)\cong 0$.

	To prove the claim for equivalences, the same argument can be adapted with the counits of the adjunctions and \eqref{eq:counit_compatibility}.
\end{proof}

\begin{proposition}
	\label{prop:ascent}
	Consider \Cref{setup:fm_cube} where $S$ is separated and $t$ is affine.
	Let $K\in D^b_{\operatorname{coh}}(Y_1\times_S Y_2)$ be relatively perfect over each $Y_i$.
	If $\Phi_K$ is fully faithful (resp.\ an equivalence) on $D^b_{\operatorname{coh}}$, then so is $\Phi_{\mathbf{L}(t^\prime)^\ast K}$.
\end{proposition}

\begin{proof}
	Denote by $K^\prime$ the kernel obtained in \Cref{prop:right_adjoint_base_change}, and by $\eta$, $\epsilon$ (resp.\ $\eta^\prime$, $\epsilon^\prime$) the unit and counit of $\Phi_K\dashv\Phi_{K^\prime}$ (resp.\ $\Phi_{\mathbf{L}(t^\prime)^\ast K}\dashv\Phi_{\mathbf{L}(t^\prime)^\ast K^\prime}$).
	Let $G_i$ be a compact generator of $D_{\operatorname{qc}}(Y_i)$.
	By \Cref{lem:equivalences_induced}, full faithfulness of $\Phi_K$ is equivalent to $\eta_{G_1}$ being an isomorphism.
	Then \eqref{eq:unit_compatibility}, together with the fact that $\lambda^K$ and $\alpha^{K^\prime}$ are isomorphisms, implies that $\eta^\prime_{\mathbf{L}t_1^\ast G_1}$ is an isomorphism.
	Since $t_1$ is affine, $\mathbf{L}t_1^\ast G_1$ is a compact generator \cite[\href{https://stacks.math.columbia.edu/tag/0BQT}{Tag 0BQT}]{StacksProject}.
	By \Cref{lem:relative_perfectness_affine_base_change}, $\mathbf{L}(t^\prime)^\ast K$ is relatively perfect over each $Y_i\times_S T$.
	As $T$ is separated (\Cref{rem:fm_cube_separated}), \Cref{lem:equivalences_induced} implies that $\Phi_{\mathbf{L}(t^\prime)^\ast K}$ is fully faithful.

	If $\Phi_K$ is an equivalence, then $\epsilon_{G_2}$ is an isomorphism.
	Then \eqref{eq:counit_compatibility} implies that $\epsilon^\prime_{\mathbf{L}t_2^\ast G_2}$ is an isomorphism.
	Since $\mathbf{L}t_2^\ast G_2$ is a compact generator, \Cref{lem:equivalences_induced} implies that $\Phi_{\mathbf{L}(t^\prime)^\ast K}$ is an equivalence.
\end{proof}

\begin{theorem}
	\label{thm:fiberwise}
	Consider \Cref{setup:fm_cube} with $S$ separated.
	Let $K\in D^b_{\operatorname{coh}}(Y_1\times_S Y_2)$ be relatively perfect over each $Y_i$.
	Then the following are equivalent:
	\begin{enumerate}
		\item \label{thm:fiberwise1} $\Phi_{K}$ is fully faithful (resp.\ an equivalence) on $D^b_{\operatorname{coh}}$
		\item \label{thm:fiberwise2} $\Phi_{\mathbf{L}(t^\prime)^\ast K}$ is fully faithful (resp.\ an equivalence) on $D^b_{\operatorname{coh}}$ for every affine morphism $t\colon T\to S$ from a Noetherian scheme
		\item \label{thm:fiberwise3} $\Phi_{\mathbf{L}(t^\prime)^\ast K}$ is fully faithful (resp.\ an equivalence) on $D^b_{\operatorname{coh}}$ for every morphism $t\colon \operatorname{Spec}(k)\to S$ from a field
		\item \label{thm:fiberwise4} $\Phi_{\mathbf{L}(t^\prime)^\ast K}$ is fully faithful (resp.\ an equivalence) on $D^b_{\operatorname{coh}}$ for every closed point $p\in S$ with associated closed immersion $t\colon \operatorname{Spec}(\kappa(p))\to S$
		\item \label{thm:fiberwise5} $\Phi_{\mathbf{L}(t^\prime)^\ast K}$ is fully faithful (resp.\ an equivalence) on $D^b_{\operatorname{coh}}$ for some affine surjective morphism $t\colon T\to S$ from a Noetherian scheme.
	\end{enumerate}
\end{theorem}

\begin{proof}
	By \Cref{prop:ascent}, $\eqref{thm:fiberwise1}\implies \eqref{thm:fiberwise2}$.
	By \cite[\href{https://stacks.math.columbia.edu/tag/01SI}{Tag 01SI}]{StacksProject}, $\eqref{thm:fiberwise2}\implies \eqref{thm:fiberwise3}$.
	It is obvious that $\eqref{thm:fiberwise3}\implies \eqref{thm:fiberwise4}$.
	By \Cref{prop:descent}\eqref{prop:descent2}, $\eqref{thm:fiberwise4} \implies \eqref{thm:fiberwise1}$.
	Taking $t=1_S$ shows $\eqref{thm:fiberwise1}\implies \eqref{thm:fiberwise5}$.
	We prove that $\eqref{thm:fiberwise5}\implies \eqref{thm:fiberwise4}$.
	Let $t\colon T\to S$ be as in \eqref{thm:fiberwise5}.
	For each closed point $p\in S$, the scheme $T\times_S \operatorname{Spec}(\kappa(p))$ is Noetherian, and the canonical morphism $T\times_S \operatorname{Spec}(\kappa(p))\to T$ is affine.
	Moreover, $T$ is separated (\Cref{rem:fm_cube_separated}).
	By \Cref{prop:ascent}, we can ascend full faithfulness or equivalences from $T$ to $T\times_S \operatorname{Spec}(\kappa(p))$.
	The canonical morphism $T\times_S \operatorname{Spec}(\kappa(p)) \to \operatorname{Spec}(\kappa(p))$ is affine and faithfully flat.
	Hence, \Cref{prop:descent}\eqref{prop:descent1}, applied to the families $Y_i\times_S \operatorname{Spec}(\kappa(p))\to \operatorname{Spec}(\kappa(p))$, allows us to descend full faithfulness or equivalences from $T\times_S \operatorname{Spec}(\kappa(p))$ to $\operatorname{Spec}(\kappa(p))$.
	This completes the proof.
\end{proof}

\begin{corollary}
	\label{cor:orlov}
	Let $Y_1$ and $Y_2$ be proper schemes over a field $k$.
	Suppose $K\in D^b_{\operatorname{coh}}(Y_1\times_k Y_2)$ is relatively perfect over each $Y_i$.
	Let $L/k$ be a field extension and $t^\prime\colon Y_1\times_k Y_2\times_k \operatorname{Spec}(L) \to Y_1\times_k Y_2$ the canonical morphism.
	Then $\Phi_K$ is fully faithful (resp.\ an equivalence) on $D^b_{\operatorname{coh}}$ if, and only if, $\Phi_{\mathbf{L}(t^\prime)^\ast K}$ is fully faithful (resp.\ an equivalence) on $D^b_{\operatorname{coh}}$.
\end{corollary}

\begin{proof}
	Each $Y_i$ is flat over $k$, the scheme $\operatorname{Spec}(k)$ is separated, and the morphism $\operatorname{Spec}(L)\to\operatorname{Spec}(k)$ is affine and faithfully flat.
	Hence, the claim follows from \Cref{prop:ascent} and \Cref{prop:descent}\eqref{prop:descent1}.
\end{proof}

\bibliographystyle{alpha}
\bibliography{mainbib}

\end{document}